\documentclass{article}
\usepackage[utf8]{inputenc}
\usepackage{graphicx} 
\usepackage[english]{babel}
\usepackage[left=3cm,top=4cm,right=3cm,bottom=4cm]{geometry}
\usepackage{amsthm}
\usepackage[tbtags]{amsmath}
\usepackage{amssymb}
\usepackage{bbm}
\usepackage{dsfont}
\usepackage{tikz}
\usepackage{enumitem}
\usepackage{mathtools}
\usepackage{stmaryrd}
\usepackage{authblk}
\usepackage{todonotes}
\usepackage{subcaption}
\usepackage{comment}
\usepackage{cases}
\usepackage{empheq}
\usepackage{hyperref}
\usepackage[labelsep=period]{caption} 
\usepackage[backend=biber,style=numeric-comp,sorting=none,maxbibnames=99,firstinits=true]{biblatex}
\allowdisplaybreaks

\newcommand{\Dphi}{D_{\varphi}}

\providecommand{\Supp}{\operatorname{supp}}                            
\providecommand{\supp}{\Supp}

\providecommand{\diam}{\operatorname{diam}}

\providecommand{\Bd}{{\boldsymbol{d}}}
\providecommand{\Be}{{\boldsymbol{e}}}

\providecommand{\Bn}{{\boldsymbol{n}}}

\providecommand{\Bx}{{\boldsymbol{x}}}
\providecommand{\By}{{\boldsymbol{y}}}

\providecommand{\BY}{{\boldsymbol{Y}}}

\newcommand{\Bdelta}     {\boldsymbol{\delta}}

\newcommand{\Beta}       {\boldsymbol{\eta}}                

\newcommand{\dd}{\,\mathrm{d}} 

\newtheorem{theorem}{Theorem}[section]
\newtheorem{corollary}{Corollary}[section]
\newtheorem{lemma}{Lemma}[section]
\newtheorem{proposition}{Proposition}[section]
\newtheorem{assumption}{Assumption}[section]
\numberwithin{equation}{section}

\theoremstyle{definition}
\newtheorem{definition}{Definition}[section]
\theoremstyle{remark}
\newtheorem{remark}{Remark}[section]

\usetikzlibrary{shapes.misc}
\tikzset{cross/.style={cross out, draw=black, minimum size=2*(#1-\pgflinewidth), inner sep=0pt, outer sep=0pt},cross/.default={2.5pt}}

\usetikzlibrary{backgrounds}

\title{Wavenumber-Explicit Well-Posedness of Bayesian Shape Inversion in Acoustic Scattering}
\author[1]{S. Kuijpers}
\author[2]{L. Scarabosio\thanks{Corresponding author: laura.scarabosio@ru.nl}}
\affil[1]{Bernoulli Institute, Faculty of Science and Engineering, University of Groningen, The Netherlands}
\affil[2]{Institute for Mathematics, Astrophysics and Particle Physics, Radboud University, Nijmegen, The Netherlands}

\date{\today}

\newcommand{\Onorm}{\left|\left(\lVert\mathcal{O}_k\rVert_{V^{\ast}}\right)_{k=1}^K\right|}

\begin{document}

\maketitle

\begin{abstract}
We consider the Bayesian approach to the inverse problem of recovering the shape of an object from measurements of its scattered acoustic field. Working in the time-harmonic setting, we focus on a Helmholtz transmission problem and then extend our results to an exterior Dirichlet (sound-soft) problem. It is well known that higher frequencies yield higher resolution but greater sensitivity to noise; here we give the first rigorous results quantifying how this sensitivity to noise depends on the wavenumber. We model the scatterer as star-shaped, with a prior on its boundary given by a series expansion of the angle-dependent radius with uniformly distributed coefficients. Our main result establishes well-posedness of the Bayesian shape inverse problem with constants explicit in the wavenumber, under problem-specific conditions on the material parameters that exclude quasi-resonant regimes. Stability estimates in the Hellinger and 1-Wasserstein metrics show that the stability constant of the posterior with respect to the data grows exponentially with the square of the wavenumber, whose magnitude must be understood not in absolute terms but relative to the spatial scale of the problem. Numerical experiments illustrate this effect.

\medskip
\noindent \textbf{Keywords}: inverse scattering problem, Bayesian inversion, shape uncertainty, Helmholtz equation, wavenumber-explicit estimates.
\end{abstract}

\section{Introduction}
%
%



Many technological applications entail acquiring information about the shape of an object from indirect measurements of its scattering behavior. Examples can be found in radar and sonar imaging \cite{borden2001mathematical,jia2002time,cervenka1993sidescan,ferguson2005application}, medical imaging \cite{adams2020uncertain,ambellan2019statistical}, and nano-optics \cite{raziman2013polarisation,sannomiya2009optical,repan2022exploiting}, to name a few. There, the uncertainty in the shape is usually quite large before acquiring any measurement. In some cases, an underlying but inaccessible ``true" shape exists, and measurements are meant to reduce the uncertainty about it; yet, due to the ill-posedness of the inversion procedure, some uncertainty will persist and one may want to quantify it to know how much the reconstructed shape can be trusted. In other cases, the shape itself is subject to intrinsic variations and one might want to infer a probabilistic description of these variations from some observations. An example of the first situation is radar imaging \cite{borden2001mathematical}, where a single object with a deterministic shape is to be retrieved. For the second situation, one could imagine having a population of objects with intrinsic shape variability, for example, due to manufacturing defects \cite{raziman2013polarisation}, and measurements of the scattering behavior of each of them. In both scenarios, the mathematical methodology that allows for quantifying the remaining shape uncertainty after data acquisition is to solve the Bayesian inverse problem for the shape. Additionally, in all applications mentioned, we know the physics underlying the measurements, namely a scattering phenomenon, so we can use it to relate the shape to the measurement data. This is the setting considered in this paper: the Bayesian shape inverse problem where the forward model describing the physics is time-harmonic wave scattering, modeled by the Helmholtz equation. We assume the measurement data to be finite-dimensional. Our main goal is to study the well-posedness of this inverse problem, and in particular to quantify a well-known qualitative phenomenon: the higher sensitivity of the inversion result to noise as the frequency increases. In the Bayesian framework, this amounts to obtaining wavenumber-explicit constants for the stability of the shape posterior with respect to the data.

The problem of reconstructing a single shape in time-harmonic wave scattering has a long tradition in applied mathematics and computational physics, explored both theoretically, e.g., studying uniqueness of the reconstruction, and numerically, through the development and analysis of reconstruction algorithms; we refer to the monograph \cite{colton1998inverse} and the more recent survey \cite{colton2018looking} for an overview and pointers to the relevant literature. In the Bayesian setting, a prior on the shape is instead updated via the data to a posterior, the two approaches being related in that a point estimate of the posterior (e.g., its mean or mode) can coincide with the solution to a deterministic inverse problem in some settings \cite{stuart2010inverse}.

To the authors' knowledge, one of the first papers addressing the Bayesian shape inverse problem is \cite{bui2014analysis}, studying the well-posedness (with no explicit tracking of the wavenumber) and numerical approximation of a sound-soft acoustic scattering problem with a lognormal prior for the boundary of a star-shaped scatterer. While \cite{bui2014analysis} considers point measurements, well-posedness with far-field data is studied in \cite{huang2021bayesian}. In \cite{carpio2020bayesian}, the authors build a prior using topological sensitivity analysis, handling a possibly uncertain number of scattering objects via model selection; outside scattering applications, more unstructured topological changes have been tackled via a Bayesian level set approach \cite{iglesias2016bayesian,dunlop2017hierarchical}. On the numerical side, Bayesian shape inversion has been studied for limited-aperture data via a procedure combining classical and Bayesian approaches \cite{li2020extended}, and for phaseless far-field data \cite{yang2020bayesian}. When the noise level is not too small \cite{schillings2016scaling}, that is, the posterior is not too concentrated, moments with respect to the posterior can be computed efficiently using high-order quadrature rules \cite{dolz2022isogeometric, henriquez2021shape}, as also done in \cite{gantner2018higher} for electrical impedance tomography. In our work, we focus on the acoustic transmission problem, sometimes also referred to as a medium scattering problem, and on a sound-soft scattering problem. We rely on an explicit parametrization for the boundary of the star-shaped scatterer to define its prior, using an expansion with uniformly distributed coefficients, similarly to \cite{gantner2018higher,henriquez2021shape} and to expansions previously suggested for uncertain coefficients \cite{schwab2012sparse}. However, while \cite{gantner2018higher, dolz2022isogeometric, henriquez2021shape} focus their analysis on a reference geometric configuration, here we work on the physical configuration, and as such do not rely on higher-order smoothness of the map from the shape to the noiseless measurements.

With the inverse problem setting at hand, we study, analytically and numerically, 
the dependence of the posterior distribution on the wavenumber. Under suitable conditions on the material parameters that exclude quasi-resonant regimes, we prove the existence and uniqueness of the posterior distribution for all frequencies, and we derive wavenumber-explicit estimates on the stability of the posterior with respect to perturbations in the data. For the latter, we use both the Hellinger and 1-Wasserstein distances, the second being potentially better suited to the multimodal posteriors that can arise at high frequencies and large noise levels. As shown in Theorems \ref{thm:wellposed_HSTP} and \ref{thm:wellposed_soft}, the posterior stability constant grows exponentially with the square of the wavenumber. We then report numerical experiments that provide further intuition on the analytical results.

Our work falls in the latest developments on wavenumber-explicit uncertainty quantification for the Helmholtz equation. The work \cite{moiola2019acoustic} has provided wavenumber-explicit estimates that are particularly suited for shape uncertainty quantification. In particular, the constants in the estimates depend only on the size of scatterer but not its shape. These results and techniques have been extended to related settings \cite{graham2019helmholtz, chaumont2023explicit} and constitute the starting point for wavenumber-explicit estimates for the forward propagation of uncertainty \cite{pembery2020helmholtz, spence2023wavenumber}, with \cite{hiptmair2025frequency} particularly addressing shape uncertainty, and for our results as well. From a more algorithmic perspective, the paper \cite{ganesh2021quasi} derives wavenumber-explicit estimates for the convergence of quasi-Monte Carlo quadrature rules for a Helmholtz equation with stochastic refractive index, and \cite{scarabosio2022deep} provides computational results on the wavenumber-dependent performance of a neural network surrogate for shape uncertainty. Recent developments can also be found in forward uncertainty quantification with integral formulations \cite{dolz2024parametric} and model order reduction \cite{bonizzoni2020least}. 

This paper is structured as follows. In Section \ref{sec:Model}, we state our model problem, namely a Helmholtz transmission problem with uncertain shape. In Section \ref{sec:Bayes}, we define the Bayesian inverse problem and recall some general well-posedness results. In Section \ref{sec:well-posedness}, we establish our novel wavenumber-explicit estimates on the well-posedness of the Bayesian shape inverse problem. These are first derived for the transmission problem and then extended to sound-soft scattering. Section \ref{sec:Results} describes a computational procedure to sample from the posterior and shows numerical results that provide some practical intuition for the effect of the wavenumber. 
We conclude with some final considerations in Section \ref{sec:conclusions}.


\section{Model Problem}\label{sec:Model}
We consider a physical setting in which a time-harmonic wave is focused on an object in free space and scatters on it. We focus on the case of a penetrable particle, namely, part of the incoming wave is reflected but, also, part of it is transmitted. In Subsection \ref{ssec:helmholtzextensions}, we will extend the results to the case of an impenetrable particle with sound-soft scattering. With $d=2,3$ being the dimension of the physical space, we denote the incoming wave as $u^i\colon \mathbb{R}^d\rightarrow \mathbb{C}$ and the scattered field as $u\colon \mathbb{R}^d\rightarrow \mathbb{C}$, so that the total field is $u^T = u^i+u$. The scattered field $u$, so also the total field, satisfies the Helmholtz equation, describing the linear propagation of pressure waves in acoustics, and, for $d=2$, the propagation of transverse electric (TE) or transverse magnetic (TM) modes in electromagnetism. 

The shape of the scatterer carries some uncertainty and our goal is to characterize and reduce the latter using the information provided by measurements of the scattering behavior. 

In Subsection \ref{ssec:shapemodeling} we discuss a probabilistic model for the shape, and in Subsection \ref{ssec:helmholtz} we outline the Helmholtz transmission problem describing the physics. The latter will be used, in Section \ref{sec:Bayes}, to relate the uncertain shape to the measurement data.

\begin{figure}[t]
\centering
          \begin{tikzpicture}[scale=0.9,framed, inner frame xsep=0.2cm, inner frame ysep=0.4cm]
        \filldraw[fill=cyan!30!white,draw=black!60!cyan,domain=0:6.28,samples=100]  plot ({1.6*cos(\x r)*(1-0.08*cos(3*\x r)+0.04*cos(6*\x r)+0.01*sin(11*\x r)}, {1.6*sin(\x r)*(1-0.08*cos(3*\x r)+0.04*cos(6*\x r)+0.01*sin(11*\x r)});
        \draw[very thin] (0,0)->(1.3,1);
        \node at (0.35,0.9) {$r(\omega,\varphi)$};
        \node at (-0.2,-0.7) {$D_{in}(\omega)$};
        \node at (2,1.6) {$D_{out}(\omega)$};
        \node at (1.5,-1.4) {$\Gamma(\omega)$};
        \draw[thick,->] (-3.8,0)--(-2.2,0);
        \draw[semithick] (-2.8,-0.2)--(-2.8,0.2);
        \draw[semithick] (-2.7,-0.2)--(-2.7,0.2);
        \node at (-3.4,0.25) {$u_i$};
    \end{tikzpicture}  
\caption{Description of the physical setting (in two space dimensions): an incoming plane wave scatters on a star-shaped object.}\label{fig:particle}
\end{figure}

\subsection{Shape Modeling}\label{ssec:shapemodeling}
We introduce the probability space $(\Omega,\mathcal{A},\mathds{P})$, where $\Omega$ is the set of all possible shape realizations of the particle, $\mathcal{A}\subseteq\mathcal{P}(\Omega)$ is a $\sigma$-algebra on $\Omega$ and $\mathds{P}$ is the probability measure on $\mathcal{A}$. For $\mathbb{P}$-a.e. $\omega\in\Omega$, we denote by $D_{in}(\omega)\subseteq \mathbb{R}^d$, $d=2,3$, the bounded domain corresponding to the scattering object and by $\Gamma(\omega)$ its boundary. The outer domain is denoted by $D_{out}(\omega)$, such that $D_{in}(\omega)\cup\Gamma(\omega)\cup D_{out}(\omega)=\mathbb{R}^d$, for every $\omega\in\Omega$.\\

We want to be able to parameterize the boundary of the scatterer while being in a setting where the wavenumber-explicit estimates from \cite{moiola2019acoustic} hold. Therefore, we make the following assumption on the scatterer's shape.
\begin{assumption}\label{ass:star-shaped}
    For $\mathbb{P}$-a.e. $\omega\in\Omega$, the bounded domain of the scattering object $D_{in}(\omega)$ is star-shaped with respect to the origin, and its boundary $\Gamma=\Gamma(\omega)$ is a measurable map from $\Omega$ to a separable subspace of closed, $C^{0,1}$-curves.
\end{assumption}
We remind that $D_{in}(\omega)$ being star-shaped with respect to the origin means that, whenever $\Bx\in D_{in}(\omega)$, the line segment $[\mathbf{0},\Bx]\subseteq D_{in}(\omega)$. 
Assumption \ref{ass:star-shaped} allows us to use an angle-dependent radius to parametrize the boundary of the scatterer, in the same spirit as in \cite{hiptmair2018large}. Namely, for $d=2$ we can write
\begin{equation}\label{eq:gamma2d}
    \Gamma(\omega) = \left\{\Bx=(r(\omega,\varphi)\cos(\varphi),r(\omega,\varphi)\sin(\varphi))\in\mathbb{R}^2, \text{ such that } \varphi\in \Dphi:=[0,2\pi)\right\},
\end{equation}
where $r=r(\omega,\varphi)$ is the radius. For $d=3$, using the notation $\varphi:=(\varphi_1,\varphi_2)$, we write
\begin{equation}\label{eq:gamma3d}
\begin{split}
    \Gamma(\omega) = &\left\{\Bx=(r(\omega,\varphi)\cos(\varphi_1)\sin(\varphi_2),r(\omega,\varphi)\sin(\varphi_1)\sin(\varphi_2), r(\omega,\varphi)\cos(\varphi_2))\in\mathbb{R}^3,\right.\\
    &\left. \;\text{ such that } \varphi\in \Dphi:=[0,2\pi)\times [0,\pi]\right\}.    
\end{split}
\end{equation}
The symbol $\Dphi$ has been introduced to treat the angle dependency in a unified way for $d=2,3$. Our notation is summarized in Figure \ref{fig:particle}. 

The setting introduced allows us to work with a function, the radius, to characterize shapes. We model it probabilistically as follows.
\begin{assumption}\label{ass:Expansion}
    For $\mathbb{P}$-a.e. $\omega\in\Omega$, the angle-dependent radius has a representation of the form
    \begin{equation}
        r(\omega,\varphi) = r_0(\varphi)+\sum^\infty_{j=1}\beta_jY_{j}(\omega)\psi_j(\varphi),\quad\varphi\in \Dphi,
        \label{eq:radiusomega}
    \end{equation}
    where $r_0\in C^{0,1}_{per}(\Dphi)$ is uniformly bounded from below by a positive constant, $\left\{\psi_j\right\}_{j\geq 1}\!\subset \!C^{0,1}_{per}(\Dphi)$ are linearly independent functions normalized to $\lVert \psi_j\rVert_{C^{0,1}_{per}(\Dphi)}=1$ for all $j\geq 1$, $\{Y_j\}_{j\geq 1}$ are random variables, and $(\beta_j)_{j\geq 1}$ is an absolutely convergent sequence of real numbers (that is, $(\beta_j)_{j\geq 1}\in \ell^1(\mathbb{N})$).
\end{assumption}
Later in this subsection, we comment on possible choices for $\left\{\psi_j\right\}_{j\geq 1}$.

The following assumption prescribes a probability distribution on the random coefficients in \eqref{eq:radiusomega} and a condition on the deterministic coefficients ensuring uniform positivity and boundedness of the radius. In Section \ref{sec:Bayes}, this will be the \textsl{prior} distribution.
\begin{assumption}\label{ass:uniform}
    The random variables $\{Y_j\}_{j\geq 1}$ are independent, identically distributed $\sim \mathcal{U}([-1, 1])$, and the deterministic coefficients are such that $\sum_{j\geq 1} |\beta_j| \leq \gamma_{\beta} r_0^{-}$, for some fixed $\gamma_{\beta}\in (0,1)$ and $r_0^{-}:=\inf_{\varphi\in \Dphi}r_0(\varphi)$.
\end{assumption}
In particular, the factor $\gamma_{\beta}$ determines how much variation around the mean shape is allowed.

Given the expansion \eqref{eq:radiusomega} for a fixed sequence $(\beta_j)_{j\geq 1}$ and random variables $\left\{Y_j\right\}_{j\geq 1}$ fulfilling Assumption \ref{ass:uniform}, we define 
\begin{equation}\label{eq:radiusspace}
    X:=\left\{r(\omega)\text{ given by }\eqref{eq:radiusomega} , \,\omega\in\Omega \right\},
\end{equation}
containing all possible realizations of the radius. 

We consider $X$ to be a measure space via the \textsl{pushforward} of the \textsl{product measure} $\bigotimes_{j\geq 1} \mu_j$  under \eqref{eq:radiusomega}, where $\mu_j=\mathcal{U}([-1,1])$, $j\geq 1$. On the other hand, for the estimates in Section \ref{sec:well-posedness}, we also consider $X$ to be a subset of the Banach space $C^{0,1}_{per}(\Dphi)$, in particular, it is a bounded set contained in the ball of center $r_0$ and radius $\sum_{j\geq 1}|\beta_j|$. Such double identification for $X$ is possible because, under Assumptions \ref{ass:Expansion}-\ref{ass:uniform}, the parametrization \eqref{eq:radiusomega} establishes a continuous map from the product topology to the $C^{0,1}_{per}(\Dphi)$-topology.

\smallskip
In the expansion \eqref{eq:radiusomega}, 
$r_0\in C^{0,1}_{per}(\Dphi)$ corresponds to a mean shape, 
and
the summation in \eqref{eq:radiusomega} is a random field with zero average. A natural choice for the latter is a Karhunen-Lo\`eve expansion, in which case $\left\{\psi_j\right\}_{j\geq 1}$ are the eigenfunctions of the covariance operator. If one assumes the covariance kernel of the random field to be rotationally invariant, which might be reasonable for many applications, then it is easy to see that $\left\{\psi_j\right\}_{j\geq 1}$ coincide with Fourier harmonics for $d=2$, see \eqref{eq:fourierpsi}, 
and spherical harmonics for $d=3$.

When the nominal radius $r_0$ is constant, using Fourier or spherical harmonics corresponds to modeling variations in the normal direction using the eigenfunctions of the Laplace-Beltrami operator on the nominal interface. When $r_0$ is not a constant, it is still possible to use eigenvalues of the Laplace-Beltrami operator to model variations in the normal direction only \cite{church2020domain}, but their expression is not known explicitly anymore and has to be approximated numerically \cite{nasikun2018fast}.

The possibilities for $\left\{\psi_j\right\}_{j\geq 1}$ discussed so far fall under the category of spectral expansions. Other modeling choices are possible, for instance, considering functions with localized supports \cite{van2023exploiting,zech2018sparse}, possibly constructed in such a way to have some prescribed statistical properties \cite{bachmayr2018representations,bachmayr2023multilevel}. The resulting prior in this case would be in the spirit of a Besov prior \cite{saksman2009discretization,dashti2011besov}. 

If the functions $\left\{\psi_j\right\}_{j\geq 1}$ are smooth, then the decay of the coefficient sequence $(\beta_j)_{j\geq 1}$ characterizes the smoothness of the random field for the shape variations \cite{hiptmair2018large}. In applications, it might be desirable to infer this decay from experimental data, in the spirit of a hierarchical statistical approach. Here, for simplicity, we assume the coefficient sequence to be given.

\begin{remark}[Radius regularity]
    Assumptions \ref{ass:Expansion} and \ref{ass:uniform} guarantee in particular a realization-independent upper bound on $\lVert r\rVert_{C^{0,1}_{per}(\Dphi)}$, $r\in X$.
\end{remark}

Some of the results in Section \ref{sec:well-posedness} hold under the assumption that all scatterer realizations are star-shaped with respect to a ball centered at the origin rather than the origin only. We remind that a domain is star-shaped with respect to a ball if it is star-shaped with respect to every point in that ball. The following lemma shows that this condition is always fulfilled in our setting.

\begin{lemma}\label{lem:starshaped}
  If Assumptions \ref{ass:Expansion}-\ref{ass:uniform} hold, then there exists $\tilde{\gamma}>0$ such that, for $\mathbb{P}$-a.e. $\omega\in\Omega$, $D_{in}(\omega)$ is star-shaped with respect to a ball of radius $\tilde{\gamma}\mathrm{diam}(D_{in}(\omega))$. 

  For $d=2$, we can take $\tilde{\gamma}=\left(\frac{(1-\gamma_{\beta})r_0^{-}}{\sqrt{2}(1+\gamma_{\beta})\lVert r_0\rVert_{C^{0,1}_{per}(\Dphi)}}\right)^2$, with $\gamma_{\beta}$ and $r_0^{-}$ as in Assumption \ref{ass:uniform}.

  For $d=3$, $\tilde{\gamma}=\frac{1}{2}\left(\frac{(1-\gamma_{\beta})r_0^{-}}{\sqrt{2}(1+\gamma_{\beta})\lVert r_0\rVert_{C^{0,1}_{per}(\Dphi)}}\right)^2 \left(\left(\frac{1}{(1-\gamma_{\beta})r_0^{-}}+1\right)(1+\gamma_{\beta})\lVert r_0\rVert_{C^{0,1}_{per}(\Dphi)}+1\right)^{-1}$.
\end{lemma}
\begin{proof}
According to Lemma 2.1 in \cite{moiola2019acoustic}, $D_{in}(\omega)$, $\omega\in\Omega$, is star-shaped with respect to a ball $B_{a}(\boldsymbol{0})$ of radius $a$ centered at the origin if and only if it is Lipschitz and $\Bx\cdot \Bn(\Bx)\geq a$ for all $\Bx\in\Gamma(\omega)$ for which the outer normal $\Bn$ to the interface is defined, see also Lemma 5.4.1 in \cite{moiola2011trefftz}. The result then follows by directly computing $\Bx\cdot\Bn(\Bx)$ for the parametrizations \eqref{eq:gamma2d} and \eqref{eq:gamma3d}, for $d=2$ and $d=3$ respectively. For $d=2$ one can use polar coordinates, for $d=3$ it is more convenient for calculations to use the implicit parametrization $x_1^2+x_2^2+x_3^2=r^2$ for the surface. In order to obtain a value $\tilde{\gamma}$ independent of $\omega\in\Omega$, it is sufficient to use the lower and upper bounds on the radius $r(\omega,\cdot)$ and its first derivatives ensured by Assumption \ref{ass:uniform}.
\end{proof}

\subsection{Helmholtz Scattering Transmission Problem}\label{ssec:helmholtz}
We now introduce the mathematical model which will be used in the statistical model in Section \ref{sec:Bayes} to connect our unknown, the shape, to the measurement data. For later convenience, we define
\begin{equation}\label{eq:holdall}
    D_{in,H}:=\cup_{\omega\in\Omega} D_{in}(\omega),\quad D_{out,H}:=\cup_{\omega\in\Omega} D_{out}(\omega),
\end{equation}
to be the hold-all domains for the scatterer realizations and the exterior domains, respectively. The following is based on Definition 2.4 in \cite{moiola2019acoustic}.

\begin{definition}\label{def:HSTP} (Scattering transmission problem)\\
    Let $\kappa_0 \in \mathbb{R}_{>0}$ be the wavenumber, and let $n_{in}$, $n_{out}$, $\alpha_{in}$, $\alpha_{out}\in \mathbb{R}_{>0}$ be material parameters. For $\omega\in\Omega$, define the piecewise constant functions
    \begin{equation}\label{eq:alphakappacoeffs}
        \alpha(\omega,\Bx) \coloneqq
        \begin{cases}
            \alpha_{in}\quad&\text{if }\Bx\in D_{in}(\omega),\\
            \alpha_{out}\quad&\text{if }\Bx\in D_{out}(\omega),
        \end{cases}\qquad n(\omega,\Bx) \coloneqq
        \begin{cases}
            n_{in}\quad&\text{if }\Bx\in D_{in}(\omega),\\
            n_{out}\quad&\text{if }\Bx\in D_{out}(\omega).
        \end{cases}
    \end{equation}
    Let $u^i$ be a solution of $\alpha_{out}\Delta u^i+\kappa_0^2n_{out}u^i=0$, that is $C^\infty$ in a neighbourhood of $\overline{D_{in,H}}$. For $\omega\in\Omega$, the scattered field $u=u(\omega)$ satisfies
    \begin{subequations}\label{eq:HSTP}
\begin{empheq}[left={\empheqlbrace}]{alignat=2}
        &-\nabla\cdot (\alpha(\omega,\cdot)\nabla u) - \kappa_0^2 n(\omega,\cdot)u =f(\omega,\cdot)\quad &&\text{in }D_{in}(\omega)\cup D_{out}(\omega),\label{eq:helmholtz}\\
        &\llbracket u \rrbracket = 0 \quad && \text{on }\Gamma(\omega),\label{eq:jump1}\\
        &\left\llbracket \alpha(\omega,\cdot)\frac{\partial u}{\partial \Bn}\right\rrbracket = - \left\llbracket \alpha(\omega,\cdot)\right\rrbracket \frac{\partial u^i}{\partial \Bn}\quad && \text{on }\Gamma(\omega),\label{eq:jump2}\\       
        &\lim_{\lVert\Bx\rVert\rightarrow\infty}\lVert\Bx\rVert^{\frac{d-1}{2}}\left(\frac{\partial}{\partial \lVert\Bx\rVert}-i\kappa_0\sqrt{n_{out}/\alpha_{out}}\right)u = 0,&&\label{eq:Sommerfeld}
\end{empheq}
\end{subequations}
    with $f(\omega,\cdot)=\nabla\cdot (\alpha(\omega,\cdot)\nabla u^i) + \kappa_0^2 n(\omega,\cdot)u^i$. Here $\llbracket\cdot\rrbracket$ denotes the jump at the boundary $\Gamma(\omega)$, $\Bn$ is the unit normal vector to $\Gamma(\omega)$ pointing into $D_{out}(\omega)$ and $\frac{\partial}{\partial \Bn}$ the corresponding Neumann trace.
\end{definition}
A notable case is a plane wave excitation, $u^i(\Bx) =\exp(i\kappa_0\sqrt{n_{out}/\alpha_{out}}\Bd\cdot\Bx)$, where $\Bd$ is a vector with unit Euclidean norm describing the direction of propagation and $i=\sqrt{-1}$.  

Equation \eqref{eq:helmholtz} is the Helmholtz equation. The two jump conditions \eqref{eq:jump1}-\eqref{eq:jump2} guarantee continuity of the Dirichlet and Neumann traces of the total field $u+u^i$ across $\Gamma(\omega)$, respectively. 
Equation \eqref{eq:Sommerfeld} is the Sommerfeld radiation condition, which guarantees that the scattered field only radiates outwards from known sources and does not allow for arbitrary waves coming in from infinity. By rescaling, some of the parameters $\alpha_{in}$, $\alpha_{out}$, $n_{in}$, $n_{out}$ in \eqref{eq:HSTP} are redundant. We keep them in, however, for physical intuition and straightforward use in applications. Although $\kappa_0$ could in principle be complex with positive imaginary part, here we confine ourselves to real wavenumbers in order to use the estimates from \cite{moiola2019acoustic}.

We will be working in the nontrapping regime, meaning we exclude the presence of resonant geometric structures \cite{moiola2014helmholtz}. 
In our setting, this will be ensured by all scatterer realizations being star-shaped and precise relationships between the material parameters \cite{moiola2014helmholtz}. 

The following result is classical, and a proof can be found in \cite{moiola2019acoustic} (see Lemma 2.2) or references therein.

\begin{theorem}\label{theorem:solvable}
    Under Assumption \ref{ass:star-shaped}, so also under Assumptions \ref{ass:Expansion}-\ref{ass:uniform}, the scattering transmission problem as defined in Definition \ref{def:HSTP} admits a unique solution $u\in H^1_{loc}(\mathbb{R}^d)$ for $\mathbb{P}$-a.e. $\omega \in \Omega$.
\end{theorem}

\section{Bayesian Inverse Problem Setting}\label{sec:Bayes}
We formalize the Bayesian shape inverse problem in acoustics (or electromagnetics) in Subsection \ref{ssec:setting}, following the general setting in \cite{stuart2010inverse}. In Subsection \ref{ssec:wellposedgeneral}, we recall some general results on its well-posedness on which we will build in Section \ref{sec:well-posedness}.

\subsection{Definition of the Bayesian Shape Inverse Problem}\label{ssec:setting}

We assume that our data follow the statistical model
\begin{equation}
    \boldsymbol{\delta} = \mathcal{G}(r) + \boldsymbol{\eta},
    \label{eq:AdditiveNoise}
\end{equation}
with finite-dimensional measurements $\Bdelta\in\mathbb{R}^K$, $K\in\mathbb{N}$, modeled by the output of a measurement operator $\mathcal{G}: X\rightarrow\mathbb{R}^K$, with $X$ as in \eqref{eq:radiusspace}, corrupted by additive Gaussian noise $\Beta\sim\mathcal{N}(\mathbf{0},\Sigma)$ for a known, positive definite covariance matrix $\Sigma\in\mathbb{R}^{K\times K}$. We assume the noise to be independent of $r$ and we denote by $\left|\cdot\right|$ the Euclidean norm on $\mathbb{R}^K$. 

The \textsl{measurement operator} $\mathcal{G}$  is decomposed as $\mathcal{G}=\mathcal{O}\circ G$, where $G:X\rightarrow V$ is the so-called \textsl{forward operator} describing the physics and taking values in a separable Banach space $V$, which we make precise below, and $\mathcal{O}: V\rightarrow \mathbb{R}^K$ is the so-called \textsl{observation operator}. We assume $\mathcal{O}$ to be linear and continuous, namely, $\mathcal{O}\in (V^{\ast})^K$, where the notation means $\mathcal{O}=(\mathcal{O}_k)_{k=1}^K$, with $\mathcal{O}_k\in V^{\ast}$ for $k=1,\ldots,K$ and $V^{\ast}$ the dual space of $V$.

Since $H^1_{loc}(\mathbb{R}^d)$ as from Theorem \ref{theorem:solvable} is not a Banach space, we instead work on a ball $B_R\subset\mathbb{R}^d$, centered at the origin with radius $R>0$ large enough to contain all scatterer realizations and to ensure $\mathcal{O}\in(H^1(B_R)^\ast)^K$. For instance, if the measurements are (smoothed) point evaluations of the scattered field at locations known to lie outside the scatterer, $B_R$ should contain all such locations. We replace the Sommerfeld radiation condition \eqref{eq:Sommerfeld} by an exact boundary condition on $\partial B_R$ via the Dirichlet-to-Neumann (DtN) map \cite[Sect. 2.6.3]{nedelec2001acoustic}, whose sign properties \cite[Thm. 4.3]{hiptmair2018large}, together with Theorem \ref{theorem:solvable}, guarantee that the truncated transmission problem has a unique solution in $H^1(B_R)$, coinciding with the restriction of the solution of \eqref{eq:HSTP} to $B_R$. We then take the range of $G$ to be $V=H^1(D)$, with $D$ either $B_R$ or a realization-independent subdomain thereof. The choice of $D$ depends on the physical setting and is specified in Subsection \ref{ssec:explbounds}.

In the Bayesian approach to the shape inverse problem, we prescribe a prior on the shape, here boiling down to a prior on the radius, and we use \eqref{eq:AdditiveNoise} to obtain a posterior on the radius via Bayes' rule. 

Our \textsl{prior measure}, denoted by $\mu_0$, is the pushforward of $\bigotimes_{j\geq 1} \mathcal{U}([-1,1])$, the distribution of $(Y_j)_{j\geq 1}$, under \eqref{eq:radiusomega} with Assumptions \ref{ass:Expansion}-\ref{ass:uniform} fulfilled. For $X$ as in \eqref{eq:radiusspace} as a subset of $C^{0,1}_{per}(\Dphi)$, $\mu_0(X)=1$. The Bayesian inverse problem then reads: 

\begin{center}
   Given the prior measure $\mu_0$ on the radius $r$, find the posterior $\mu^{\Bdelta}$, namely the measure on $r$ conditioned on the observations \eqref{eq:AdditiveNoise}. 
\end{center}

\subsection{General Well-Posedness Results}\label{ssec:wellposedgeneral}

Under mild assumptions on the operator $\mathcal{G}$, the Bayesian inverse problem is often naturally well-posed in the sense of Hadamard \cite{stuart2010inverse}. When the prior is defined in terms of uniform random variables, existence and uniqueness of the posterior measure are ensured by Proposition \ref{prop:existence}, while stability to data perturbations is summarized in Proposition \ref{prop:stability}.

\begin{proposition}[Theorem 6.31 in \cite{stuart2010inverse}]\label{prop:existence}
    Assume that $\mathcal{G}:X\rightarrow\mathbb{R}^K$ in \eqref{eq:AdditiveNoise} is $\mu_0$-measurable. Then the posterior measure $\mu^{\boldsymbol{\delta}}(dr)=\mathds{P}(dr|\boldsymbol{\delta})$ given data $\boldsymbol{\delta}$ is absolutely continuous with respect to the prior measure $\mu_0(dr)$ and has Radon-Nikodym derivative
    \begin{equation}
        \frac{d\mu^{\boldsymbol{\delta}}}{d\mu_0}(r) \propto \exp{(-\Phi(r;\boldsymbol{\delta}))},
        \label{eq:RadonNikodym}
    \end{equation}
    where
    \begin{equation}\label{eq:potential}
        \Phi(r;\boldsymbol{\delta}) = \frac{1}{2}|\boldsymbol{\delta}-\mathcal{G}(r)|^2_{\Sigma}
    \end{equation}
    and $|\cdot|_{\Sigma} = |\Sigma^{-\frac{1}{2}}\cdot|$ with $|\cdot|$ the Euclidean norm. 
    In particular, 
     \eqref{eq:RadonNikodym} is well-defined if $\mathcal{G}:X\rightarrow \mathbb{R}^K$ is continuous and $\mu_0(X) = 1$.
\end{proposition}

The result above justifies the application of Bayes' rule for measures on $X$:
\begin{equation}\label{eq:posterior}
    \mu^{\Bdelta}(dr)\propto \frac{d\mu^{\boldsymbol{\delta}}}{d\mu_0}(r)\mu_0(dr)\propto \exp{\left(-\frac{1}{2}|\boldsymbol{\delta}-\mathcal{G}(r)|^2_{\Sigma}\right)}\mu_0(dr),
\end{equation}
where the first proportionality constant, which is also a normalization constant, is the so-called model evidence. The Radon-Nikodym derivative \eqref{eq:RadonNikodym} corresponds to the \textsl{likelihood}, and the quantity \eqref{eq:potential} is often referred to as \textsl{potential}.

For the stability of the posterior, in the following we ask for local boundedness and local Lipschitz continuity of the potential with respect to the data. We use the notation $\lVert F \rVert_{L^p_{\nu}(W_1,W_2)}=\left(\int_{W_1} \lVert F(v)\rVert_{W_2}^p \dd\nu(v)\right)^{\frac{1}{p}}$, for $1\leq p<\infty$, $W_1$ and $W_2$ Banach spaces (or possibly subsets of them) with $W_1$ separable, and integration with respect to the measure $\nu$ on $W_1$. When $p=\infty$, $\lVert F \rVert_{L^{\infty}_{\nu}(W_1,W_2)}=\text{esssup}_{v\in W_1} \lVert F(v)\rVert_{W_2}$. When $W_2=\mathbb{R}$, we simply write $L^p_{\nu}(W_1)$.

\begin{assumption}[Assumption 2 in \cite{hoang2012sparse}]\label{ass:stability}
    Let $\mu_0$ be the prior measure on $X$. The function $\Phi:X\times\mathbb{R}^K\rightarrow\mathbb{R}$ satisfies:
    \begin{enumerate}[label=(\roman*)]
        \item for each $\gamma>0$, there is a constant $M_{\gamma}>0$ and a set $X_\gamma \subseteq X$ of positive $\mu_0$-measure such that, for all $r \in X_{\gamma}$ and for all $\Bdelta$ such that $|\Bdelta| \leq \gamma$,
        \begin{equation*}
            0 \leq \Phi(r;\boldsymbol{\delta}) \leq M_{\gamma} ;
        \end{equation*}
        \item there is a mapping $\mathcal{K}: \mathbb{R}\times X \rightarrow \mathbb{R}$ such that, for each $\gamma>0$, $\mathcal{K}(\gamma,\cdot)\in L^2_{\mu_0}(X)$, and, for every $|\boldsymbol{\delta}|,|\boldsymbol{\delta}'|\leq \gamma$, it holds
        \begin{equation*}
            |\Phi(r;\boldsymbol{\delta})-\Phi(r;\boldsymbol{\delta}')|\leq \mathcal{K}(\gamma,r)|\boldsymbol{\delta}-\boldsymbol{\delta}'|.
        \end{equation*}
    \end{enumerate}
\end{assumption}

We also need a distance between measures. Here we consider the Hellinger distance and the 1-Wasserstein distance. Given two measures $\mu,\mu'$ on the same measure space, both absolutely continuous with respect to a reference measure $\nu$, their Hellinger distance defined as
\begin{equation}
    d_{\text{Hell}}(\mu,\mu')=\sqrt{\left(\frac{1}{2}\int\left(\sqrt{\dfrac{\dd\mu}{\dd\nu}}-\sqrt{\dfrac{\dd\mu'}{\dd\nu}}\right)^2\dd\nu\right)}.
\end{equation}
As noted for instance in \cite[Sect. 6.7]{stuart2010inverse}, the Hellinger distance is useful to bound expectations and covariances of random variables under different measures.

In inverse scattering, the dependence of the scattered field on the uncertain radius $r$ is oscillatory, and the higher the frequency, the more severe the oscillations are \cite{hiptmair2025frequency}. For the inverse problem, with the forward model entering the likelihood, this means that, at higher frequencies and with noise levels which are not very small, the posterior is more prone to be multimodal. We will see the emergence of this effect in Subsection \ref{sssec:adjustedSNRanddiscr} already for frequencies that are not excessivily high. It might thus be interesting to compare the posteriors in terms of their shift with respect to each other rather than their moments. For this, we also consider stability in a Wasserstein metric, focusing on the 1-Wasserstein distance. 

The 1-Wasserstein distance between two probability measures $\mu$ and $\mu'$ on a metric space $(X, d_X)$ is 
\begin{equation}
	d_{W^1}(\mu,\mu') = \inf_{\pi\in \Pi(\mu,\mu')}\int_X d_X(r,r')\pi(\mathrm{d} r \mathrm{d}r'),
\end{equation}
with $\Pi(\mu,\mu')$ the set of all couplings of $\mu$ and $\mu'$. In the following, we denote $\lVert\mu_0\rVert_{\mathcal{P}^2}^2\!=\!\inf_{r\in X}\!\int_X d_X^2(r, r')\mu_0(\mathrm{d} r')$.

\begin{proposition}[Proposition 3 in \cite{hoang2012sparse} and Corollary 19 in \cite{sprungk2020local}]\label{prop:stability}
    Under Assumption \ref{ass:stability}, $\mu^{\boldsymbol{\delta}}$ in \eqref{eq:posterior} is locally Lipschitz continuous on the data $\boldsymbol{\delta}$ with respect to the Hellinger and 1-Wasserstein distances: for each $\gamma>0$ there are $C_{\gamma,\mathrm{Hell}}, C_{\gamma,W^1}>0$ such that, if $|\boldsymbol{\delta}|,|\boldsymbol{\delta}'|\leq \gamma$, then
    \begin{equation}\label{eq:stability}
    		d_{\mathrm{Hell}}(\mu^{\boldsymbol{\delta}},\mu^{\boldsymbol{\delta'}})\leq C_{\gamma,\mathrm{Hell}}|\boldsymbol{\delta} - \boldsymbol{\delta}'|,\qquad
    		d_{W^1}(\mu^{\boldsymbol{\delta}},\mu^{\boldsymbol{\delta'}})\leq C_{\gamma,W^1}|\boldsymbol{\delta} - \boldsymbol{\delta}'|.
    \end{equation}
For $\mathcal{K}(\gamma,\cdot)$ and $M_{\gamma}$ as in Assumption \ref{ass:stability} and constants $\Tilde{C}_{\mathrm{Hell}}, \Tilde{C}_{W^1}$ independent of $\mathcal{K}$, $M$ and $\gamma$,  $C_{\gamma,\mathrm{Hell}}= \Tilde{C}_{\mathrm{Hell}}\lVert \mathcal{K}(\gamma,\cdot)\rVert_{L^2_{\mu_0}(X)}\exp\left(M_{\gamma}\right)$ and $C_{\gamma,W^1}= \Tilde{C}_{W^1}\lVert\mu_0\rVert_{\mathcal{P}^2}\lVert \mathcal{K}(\gamma,\cdot)\rVert_{L^2_{\mu_0}(X)}\exp\left(2 M_{\gamma}\right)$.
\end{proposition}

The following result shows how the conditions in Assumption \ref{ass:stability} follow from the properties of the forward operator. 

\begin{corollary}\label{cor:pdewellposedness}
   Let the operator $\mathcal{G}: X\rightarrow \mathbb{R}^K$ be given by $\mathcal{G}=\mathcal{O} \circ G$, where $G:X\rightarrow V$ is $\mu_0$-continuous and belongs to $L^{\infty}_{\mu_0}(X,V)$, and $\mathcal{O}=(\mathcal{O}_k)_{k=1}^K\in (V^{\ast})^K$. Moreover, for additive Gaussian noise such that $\Phi(r,\boldsymbol{\delta})$ is as in \eqref{eq:potential}, let $\lambda_{min}>0$ denote the minimum eigenvalue of the covariance matrix $\Sigma$.  Then: 
   \begin{enumerate}
       \item[(i)] the posterior measure is absolutely continuous with respect to the prior and has Radon-Nikodym derivative \eqref{eq:RadonNikodym};
       \item[(ii)] Assumption \ref{ass:stability} is fulfilled with $M_{\gamma}=\frac{1}{2}C_{\gamma,G}^2$ and $\mathcal{K}(\gamma,r)=C_{\gamma,G}$ (independent of $r$), where
       \begin{equation}\label{eq:CgammaG}
        C_{\gamma,G}:=\lambda_{min}^{-1}\left(\gamma + \left|\left(\lVert\mathcal{O}_k\rVert_{V^{\ast}}\right)_{k=1}^K\right|\lVert G\rVert_{L_{\mu_0}^{\infty}(X,V)}\right),
       \end{equation}
       and \eqref{eq:stability} holds with $C_{\gamma,\mathrm{Hell}}\!=\!\Tilde{C}_{\mathrm{Hell}}C_{\gamma,G}\!\exp(\frac{1}{2}C_{\gamma,G}^2)$ and $C_{\gamma,W^1}\!=\!\Tilde{C}_{W^1}\!\lVert\mu_0\rVert_{\mathcal{P}^2}C_{\gamma,G}\!\exp(C_{\gamma,G}^2)$, for $\Tilde{C}_{\mathrm{Hell}}, \Tilde{C}_{W^1}>0$.
   \end{enumerate}
\end{corollary}
\begin{proof}
Claim $(i)$ and \eqref{eq:CgammaG} are a restatement of Proposition 8 in \cite{hoang2012sparse} when explicitly tracking the dependence of the constants on the operator norms. The expression for $M_{\gamma}$ follows from
	\begin{align*}
		\Phi(r,\Bdelta)&=\frac{1}{2}\left|\Bdelta-\mathcal{G}(r)\right|_{\Sigma}^2\leq \frac{1}{2}\left(\left|\Bdelta\right|_{\Sigma} + \left|\mathcal{G}(r)\right|_{\Sigma}\right)^2\leq \frac{1}{2}\lambda_{min}^{-2}\left(\gamma + \left|\mathcal{G}(r)\right|\right)^2,
	\end{align*}
	and $|\mathcal{G}(r)|\leq \left|\left(\lVert\mathcal{O}_k\rVert\right)_{k=1}^K\right|\lVert G\rVert_{L_{\mu_0}^{\infty}(X,V)}$. The estimate for $\mathcal{K}(\gamma,r)$ is obtained analogously once one bounds, as in \cite{hoang2012sparse}, 
	\begin{equation*}
		\left|\Phi(r,\Bdelta)-\Phi(r,\Bdelta')\right|=
		\frac{1}{2}\left|\langle \Bdelta+\Bdelta'-2\mathcal{G}(r),\Bdelta-\Bdelta'\rangle_{\Sigma}\right|\leq \lambda_{min}^{-1}\left(\gamma + |\mathcal{G}(r)|\right)|\Bdelta-\Bdelta'|.
	\end{equation*}
\end{proof}


\begin{remark}[Alternative bound]
Since, in the setting of Section \ref{sec:Model}, the metric space $X$ with distance induced by the $C^{0,1}$-norm is bounded, we could directly deduce stability in the 1-Wasserstein metric from $d_{W^1}(\mu,\mu')\lesssim d_{\mathrm{Hell}}(\mu,\mu')$, where the constant is twice the diameter of $X$ (i.e., $4\sum_{j\geq 1} |\beta_j|$). This implies a bound on $d_{W^1}(\mu,\mu')$ with the factor $\exp(M_{\gamma})$ instead of $\exp(2M_{\gamma})$ and thus a slightly more moderate growth of the constant as $M_{\gamma}$ increases. In the statements above, we kept the more general setting such that it could be easily extended to related problems with possibly unbounded $X$, e.g.,  with lognormal priors.
\end{remark}

\section{Wavenumber-Explicit Well-Posedness}\label{sec:well-posedness}
In this section, we show that the assumptions in Corollary \ref{cor:pdewellposedness} are fulfilled, under some conditions, when the forward operator maps the scatterer's boundary to the solution to the transmission problem \eqref{eq:HSTP}. This will imply that the shape inverse problem as stated at the end of Subsection \ref{ssec:setting} is well-posed, when the observation operator is continuous. When showing well-posedness, our emphasis is on making the constants, especially the one in the stability estimates \eqref{eq:stability}, explicit on the wavenumber.

In Subsection \ref{ssec:helmholtzbounds} we recall wavenumber-explicit bounds for a fixed scatterer as from \cite{moiola2019acoustic}, and in Subsection \ref{ssec:explbounds} we use these to  establish the well-posedness of the Bayesian shape inverse problem, with constants explicit in the wavenumber. In Subsection \ref{ssec:helmholtzextensions}, we extend the obtained results to sound-soft scattering.

\subsection{Preliminaries on Wavenumber-Explicit Bounds}\label{ssec:helmholtzbounds}

For later use, recall the estimates for a slightly more general scattering problem than \eqref{eq:HSTP}, and we already write them in probabilistic terms. We remind the notation \eqref{eq:holdall} and introduce the problem: for a.e. $\omega\in\Omega$, find $v=v(\omega)\in H^1(\mathbb{R}^d\setminus\Gamma(\omega))$ such that
\begin{subequations}\label{eq:HSTP2}
    \begin{empheq}[left={\empheqlbrace}]{alignat=2}
    &-\nabla\cdot (\alpha(\omega,\cdot)\nabla v) - \kappa_0^2 n(\omega,\cdot)v = f(\omega,\cdot)\quad &&\text{in }D_{in}(\omega)\cup D_{out}(\omega),\\
        &\llbracket v \rrbracket = g_D(\omega,\cdot) &&\text{on }\Gamma(\omega),\\
        &\left\llbracket \alpha(\omega,\cdot)\frac{\partial v}{\partial \Bn}\right\rrbracket = g_N(\omega,\cdot)&& \text{on }\Gamma(\omega),\\     
        &\lim_{\lVert\Bx\rVert\rightarrow\infty}\lVert\Bx\rVert^{\frac{d-1}{2}}\left(\frac{\partial}{\partial \lVert\Bx\rVert}-i\kappa_0\sqrt{n_{out}/\alpha_{out}}\right)v= 0&&,\label{eq:HSTP2_sommerfeld}
    \end{empheq}
    \end{subequations}
    where 
    \begin{equation*}
        f(\omega,\cdot)=\begin{cases}
            f_{in}(\omega) &\text{in }D_{in}(\omega),\\
            f_{out}(\omega) &\text{in } D_{out}(\omega),
        \end{cases}
    \end{equation*}
    for measurable $f_{in}:\Omega\rightarrow L^2(D_{in,H})$ and $f_{out}:\Omega\rightarrow L^2(D_{out,H})$, such that $\cup_{\omega\in\Omega}\supp f_{out}(\omega)$ is compact. We assume that $g_D$ and $g_N$ are, respectively, Dirichlet and Neumann traces of measurable random variables taking values in $H^{3/2}(D_{in,H})$, such that $g_D\in H^1(\Gamma(\omega))$ and $g_N\in L^2(\Gamma(\omega))$ for a.e. $\omega\in\Omega$. It is easy to see that problem \eqref{eq:HSTP} is a special case of \eqref{eq:HSTP2}. 


We introduce a weighted, piecewise $H^1$-norm, equivalent to the standard piecewise $H^1$-norm. For $D_1, D_2$ two disjoint subdomains of $\mathbb{R}^d$ where the material coefficients are $\alpha_{in}, n_{in}$ and $\alpha_{out}, n_{out}$ respectively, we denote the weighted norm by 
\begin{equation}\label{eq:weightednorm}
    \lVert v \rVert_{H^1_{\kappa_0,\alpha,n}(D_1\cup D_2)}^2:=  \alpha_{in}\lVert \nabla v\rVert_{L^2(D_1)}^2 + \kappa_0^2n_{in}\lVert v\rVert^2_{L^2(D_1)} + \alpha_{out}\lVert \nabla v\rVert_{L^2(D_2)}^2 + \kappa_0^2n_{out}\lVert v\rVert^2_{L^2(D_2)}.
\end{equation}

\begin{theorem}[Theorem 3.1 in \cite{moiola2019acoustic}]\label{thm:MoiolaSpence}
For $\omega\in\Omega$, let $D_{in}=D_{in}(\omega)$ be as in Assumption \ref{ass:star-shaped}. Furthermore, assume that
    \begin{equation}\label{eq:nontrappingcond}
        \frac{n_{in}}{n_{out}}\leq \frac{\alpha_{in}}{\alpha_{out}}.
    \end{equation}
    Then, for every $R>0$ such that $\cup_{\omega\in\Omega}\supp f_{out}(\omega)\subseteq B_R$ (where $B_R$ is a ball of radius $R$ centered at the origin), and denoting $D_R:=D_{out}\cap B_R$, the solution to \eqref{eq:HSTP2} when $g_D=g_N=0$ satisfies the following bound, for a.e. $\omega\in\Omega$:
    \begin{equation}\label{eq:normbound}   
    \begin{split}
         \lVert v \rVert_{H^1_{\kappa_0,\alpha,n}(D_{in}\cup D_R)}^2
     &\leq \left[\frac{4(\kappa_0 \diam(D_{in}(\omega)))^2}{\alpha_{in}} + \frac{(\kappa_0 R)^2}{n_{in}}\left(2\sqrt{\frac{n_{out}}{\alpha_{out}}} + \frac{d-1}{\kappa_0 R}\right)^2\right]\kappa_0^{-2}\lVert f_{in}\rVert^2_{L^2(D_{in}(\omega))}\\
     & + R^2 \left[\frac{4}{\alpha_{out}}+\frac{1}{n_{out}}\left(2\sqrt{\frac{n_{out}}{\alpha_{out}}} + \frac{d-1}{\kappa_0 R}\right)^2\right]\lVert f_{out}\rVert^2_{L^2(D_{out}(\omega))}.
    \end{split}
\end{equation}
\end{theorem}

\begin{theorem}[Theorem 3.2 in \cite{moiola2019acoustic}]\label{thm:MoiolaSpence2}
Additionally to Assumption \ref{ass:star-shaped}, assume that, for a.e. $\omega\in\Omega$, $D_{in}(\omega)$ is star-shaped with respect to a ball centered in the origin and with radius $\hat{\gamma}\diam(D_{in}(\omega))$ for $0<\hat{\gamma}\leq\frac{1}{2}$. Furthermore, assume that
    \begin{equation}\label{eq:nontrappingcond2}
        \frac{n_{in}}{n_{out}}<1<\frac{\alpha_{in}}{\alpha_{out}}.
    \end{equation}  
Then, for a.e. $\omega\in\Omega$, the solution to \eqref{eq:HSTP2} satisfies
    \begin{equation}\label{eq:normbound2}   
    \begin{split}
         &\lVert v \rVert_{H^1_{\kappa_0,\alpha,n}(D_{in}\cup D_R)}^2\\
     &\leq \left[\frac{4(\kappa_0 \diam(D_{in}(\omega)))^2}{\alpha_{in}} + \frac{(\kappa_0 R)^2}{n_{in}}\left(2\sqrt{\frac{n_{out}}{\alpha_{out}}} + \frac{d-1}{\kappa_0 R}\right)^2\right]\kappa_0^{-2}\lVert f_{in}\rVert^2_{L^2(D_{in}(\omega))}\\
     & + R^2 \left[\frac{4}{\alpha_{out}}+\frac{1}{n_{out}}\left(2\sqrt{\frac{n_{out}}{\alpha_{out}}} + \frac{d-1}{\kappa_0 R}\right)^2\right]\lVert f_{out}\rVert^2_{L^2(D_{out}(\omega))}\\
     & + 2\left[\frac{\diam(D_i(\omega))\alpha_{out}((3+2\hat{\gamma})\alpha_{in}+2\alpha_{out})}{\hat{\gamma}(\alpha_{in}-\alpha_{out})}\right]\lVert \nabla_T g_D\rVert_{L^2(\Gamma(\omega))}^2\\
     & + 2\left[\frac{2\kappa_0^2\diam(D_{in}(\omega))n_{out}^2}{\hat{\gamma}(n_{out}-n_{in})} + \frac{(3+\hat{\gamma})\alpha_{in}\left(n_{out}(\kappa_0 R)^2+\alpha_{out}\frac{(d-1)^2}{4}   \right)}{\hat{\gamma}\diam(D_{in}(\omega))(\alpha_{in}-\alpha_{out})}\right]\lVert g_D\rVert_{L^2(\Gamma(\omega))}^2\\
     &+ \frac{2}{\hat{\gamma}\alpha_{out}}\left[\frac{\diam(D_{in}(\omega))(4\alpha_{in}+2\alpha_{out})}{\alpha_{in}-\alpha_{out}}+\frac{2\left(n_{out}(\kappa_0 R)^2+\alpha_{out}\frac{(d-1)^2}{4}\right)}{\kappa_0^{2}\diam(D_{in}(\omega))(n_{out}-n_{in})}\right]\lVert g_N\rVert_{L^2(\Gamma(\omega))}^2,
    \end{split}
\end{equation}
where $\nabla_T$ denotes the tangential gradient and $R, D_R$ are as in Theorem \ref{thm:MoiolaSpence}.
\end{theorem}

The second result accounts for possible non-zero jumps at the interface under slightly stronger assumptions, namely the scatterer being star-shaped with respect to a ball and strict inequality in \eqref{eq:nontrappingcond2}. The first condition is always fulfilled under Assumptions \ref{ass:Expansion}-\ref{ass:uniform}, as from Lemma \ref{lem:starshaped}. 


Since the constants on the right-hand sides of \eqref{eq:normbound} and \eqref{eq:normbound2} do not depend on the shape of the scatterer, the estimates above are practical for usage in shape uncertainty quantification. 

Theorem \ref{thm:MoiolaSpence2} would allows us to obtain a wavenumber-explicit estimate for the solution $u$ to Definition \ref{def:HSTP}, but with a suboptimal dependence of the bound on the wavenumber \cite{chaumont2023explicit, moiolaspencenote}. A wavenumber-optimal estimate is obtained instead by applying Theorem \ref{thm:MoiolaSpence} to a weighted sum of the scattered and incident field, as it is done in \cite{chaumont2023explicit} for Maxwell's equations \cite{moiolaspencenote}.  

\begin{corollary}\label{cor:uscattbound}
     Let Assumption \ref{ass:Expansion} and \eqref{eq:nontrappingcond} hold. Moreover, let $R,D_R$ be as in Theorem \ref{thm:MoiolaSpence}, $R_{scatt}>\frac{\diam D_{in,H}}{2}$ and $\frac{n_{in}}{n_{out}}\leq \frac{\alpha_{in}}{\alpha_{out}}$. Then, for $\mathbb{P}$-a.e. $\omega\in\Omega$, the solution $u$ to \eqref{eq:HSTP} satisfies
    \begin{equation*}
        \lVert u \rVert_{H^1_{\kappa_0,\alpha,n}(D_{in}\cup D_R)}\!\leq \!C_{\kappa_0} C_{1}\alpha_{out}\lVert \nabla u^i\rVert_{L^2(D_{out}(\omega))} + (C_{\kappa_0} C_{2}\alpha_{out} + \sqrt{\alpha_{out}} C_{1})\lVert u^i\rVert_{L^2(D_{out}(\omega))} + \lVert u^i \rVert_{H^1_{\kappa_0,\alpha,n}},
    \end{equation*}
    where $C_{\kappa_0} = R \left[\frac{4}{\alpha_{out}}+\frac{1}{n_{out}}\left(2\sqrt{\frac{n_{out}}{\alpha_{out}}} + \frac{d-1}{\kappa_0 R}\right)^2\right]^{\frac{1}{2}}$, 
    $C_1=\frac{3}{2}\frac{1}{R_{scatt}(R-R_{scatt})}$ and $C_2=\frac{6}{(R-R_{scatt})^2}+\frac{3}{2}\frac{d-1}{R_{scatt}(R-R_{scatt})}$.
\end{corollary}

More precisely, the estimate above follows from an adaptation of the proof of Corollary 1.6 in \cite{chaumont2023explicit}. Differently from there, though, in the weighted sum $\tilde{u}=u + \chi u^i$ we used the mollifier $\chi(\Bx)=\max\left\{0,\min\left\{1,2\left(\frac{|\Bx|-R_{scatt}}{R-R_{scatt}}\right)^3-3\left(\frac{|\Bx|-R_{scatt}}{R-R_{scatt}}\right)^2+1\right\}\right\}$ instead of a piecewise linear one, and we considered the interface to be $\Gamma(\omega)$. The reason is that the formulation used for Maxwell's equations would correspond to a mixed formulation for Helmholtz, which admits milder smoothness requirements for the mollifier. The constant $C_{\kappa_0}$ in the corollary above corresponds to the square root of the second constant in \eqref{eq:normbound}. The constants $C_1$ and $C_2$ are instead upper bounds on $\lVert \nabla\chi\rVert_{L^{\infty}(D_{in}\cup D_R)}$ and $\lVert \Delta \chi\rVert_{L^{\infty}(D_{in}\cup D_R)}$, respectively.

\subsection{Well-Posedness of Bayesian Shape Inversion}\label{ssec:explbounds}

We can now show the continuity and boundedness of the forward operator in the shape inverse problem. From Corollary \ref{cor:pdewellposedness}, well-posedness of the Bayesian shape inverse problem will follow. 

We first address continuity, from which
part $(i)$ in Corollary \ref{cor:pdewellposedness} follows. For this, we prove two statements. The first one, Proposition \ref{prop:helmholtz_continuity1}, holds for the general scattering problem as in Definition \ref{def:HSTP}, but requires a smoother boundary for the scatterer compared to Theorems \ref{thm:MoiolaSpence}-\ref{thm:MoiolaSpence2}. The second one, Proposition \ref{prop:helmholtz_continuity2}, holds only when $\alpha_{in}=\alpha_{out}$ in Definition \ref{def:HSTP}, but, on the other hand, allows to accommodate star-shaped scatterers with Lipschitz boundary and a broader class of observation operators. The condition $\alpha_{in}=\alpha_{out}$ holds in many applications, for instance when the Helmholtz equation models transverse magnetic scattering and the scatterer is nonmagnetic, see \cite[Remark 2.1]{moiola2019acoustic} and \cite[p.8]{maier2007plasmonics}.



To obtain the first result on continuity, we make use of shape calculus. A similar strategy was used in \cite{bui2014analysis} for an obstacle problem with a lognormal prior on the radius. Since the presence of the (realization-dependent) interface causes a loss of global smoothness for the shape derivative, as explained for instance in \cite[Sect. 3]{harbrecht2013first} and as it will become clear from the proof below, 
we consider the case that the range of the operator $G$ is $V=H^1(B_R\setminus U)$, where $B_R$ has a radius sufficiently large to contain all realizations of the scatterer in its interior, and $U$ is the tube region 
\begin{equation}\label{eq:u}
 U:=\left\{\Bx\in B_R: r^{-}(\varphi_{\Bx})\leq \left|\Bx\right|\leq r^{+}(\varphi_{\Bx})\right\}
\end{equation}
containing all realizations of the scatterer's boundary, where $r^{-}(\varphi_{\Bx})=\inf_{\omega\in\Omega}r(\omega,\varphi_{\Bx})$, $r^{+}(\varphi_{\Bx})=\sup_{\omega\in\Omega}r(\omega,\varphi_{\Bx})$, and $\varphi_{\Bx}$ is the angle (or angles) of the point $\Bx$ in the polar (for $d=2$) or spherical (for $d=3$) coordinate system. Under Assumption \ref{ass:uniform}, the width of $U$ is determined by how much shape variation is allowed in the prior, encoded in the quantity $\gamma_{\beta}$ there. 

\begin{proposition}\label{prop:helmholtz_continuity1}
    Under Assumptions \ref{ass:Expansion}-\ref{ass:uniform}, assume that the radius $r$ in \eqref{eq:radiusomega} is $\mu_0$-a.s. of class $C^{2,1}$ and that its $C^{2,1}$-norm is uniformly bounded on $X$, i.e., $\mu_0(X\cap C^{2,1}_{per}(\Dphi))=1$ and $\operatorname{ess} \sup_{r\in X}\lVert r \rVert_{C^{2,1}_{per}(\Dphi)}<\infty$.
    Moreover, assume that the material parameters in \eqref{eq:alphakappacoeffs} fulfill \eqref{eq:nontrappingcond2}.
    
    Then the forward operator $G:X\rightarrow V$ mapping the radius to the solution to \eqref{eq:HSTP} restricted to $B_R\setminus U$, with $X$ as in \eqref{eq:radiusspace} and $V=H^1(B_R\setminus U)$, is locally Lipschitz $\mu_0$-continuous. 
\end{proposition}
\begin{proof}
In Step 1, we show that the G\^ateau differential $\delta u(r; \delta r)\in H^1(B_R\setminus \Gamma(\omega))$ for the solution to \eqref{eq:HSTP} exists, for $\mu_0$-a.e. $r\in X$ and all $\delta r \in C^{2,1}_{per}(\Dphi)$ with $\lVert \delta r\rVert_{C^{2,1}_{per}(\Dphi)}$ sufficiently small (note that, here, the space $H^1(B_R\setminus \Gamma(\omega))$ depends on $r$, namely, the location of the interface). In Step 2, we show that the map $\delta r\mapsto \delta u(r; \delta r)$, from $C^{2,1}_{per}(\Dphi)$ to $H^1(B_R\setminus \Gamma(\omega))$, is linear and continuous. The results of Step 1 and Step 2 allow us, in Step 3, to conclude the proof by an application of the mean value theorem \cite[Thm. 3.2.7]{drabek2007methods}.

\smallskip
\underline{\textsl{Step 1}} We show existence and uniqueness of $\delta u(r; \delta r)\in H^1(B_R\setminus \Gamma(\omega))$, for $\mu_0$-a.e. $r\in X$ and all $\delta r \in C^{2,1}_{per}(\Dphi)$ with sufficiently small norm. 
For $\mu_0$-a.e. $r\in X$, $r$ is $C^{2,1}$-regular by assumption, so that the solution $u$ to \eqref{eq:HSTP} is $\mu_0$-a.s. $H^2$-regular (in fact $H^3$-regular) in $D_{in}$ and $D_{out}$, see for instance \cite[Ch. 4]{mclean2000strongly}. 
This allows us to apply Theorem 4.3 in \cite{hiptmair2018shape}, characterizing the shape derivative for the transmission problem. Indeed, taking $\lVert\delta r\Vert_{C^{2,1}_{per}(\Dphi)}\!\!\leq\!\!\gamma_r \!\inf_{\omega\in\Omega,\varphi\in\Dphi}\!\! r(\omega,\varphi)$ for some $\gamma_r<1$, the assumptions of Corollary 3.1 in \cite{hiptmair2018shape} stating the smoothness of the shape derivative are fulfilled (where $v$ in the statement of that corollary corresponds to $\delta r \Be_{\rho}$ here and $\Be_{\rho}\in\mathbb{R}^d$ is the unit vector in the radial direction), guaranteeing that the assumptions of Theorem 4.3 in \cite{hiptmair2018shape} are also fulfilled. We conclude that, for $\mu_0$-a.e. $r\in X$ and $\delta r\in C^{2,1}_{per}(\Dphi)$ sufficiently small as stated before, the shape derivative $\delta u(r; \delta r)$ of the solution $u$ to \eqref{eq:HSTP} exists and it is the unique weak solution in $H^1(B_R\setminus \Gamma(\omega))$ to
    \begin{equation}\label{eq:shapeder_pde}
        \begin{cases}
        -\nabla\cdot (\alpha(\omega,\mathbf{x})\nabla \delta u) - \kappa_0^2 n(\omega,\mathbf{x})\delta u=0&\text{in }D_{in}(\omega)\cup D_R(\omega),\\
        \llbracket \delta u \rrbracket = -\delta r \Be_{\rho}\cdot \Bn\left\llbracket \frac{\partial u}{\partial \Bn}\right\rrbracket&\text{on }\Gamma(\omega),\\
        \left\llbracket \alpha(\omega,\mathbf{x})\frac{\partial \delta u}{\partial  \Bn}\right\rrbracket = \nabla_{\Gamma} \cdot \left(\llbracket \alpha \rrbracket \delta r(\Be_{\rho}\cdot \Bn)\nabla_{\Gamma} (u+u^i) \right)+(\delta r\Be_{\rho}\cdot \Bn)\kappa_0^2\llbracket n \rrbracket (u+u^i)  &\text{on }\Gamma(\omega),\\
         \alpha_{out}\dfrac{\partial \delta u}{\partial \Bn}=\operatorname{DtN}(\delta u) &\text{on }\partial B_R,
    \end{cases}
    \end{equation}
    see also \cite[Table 6]{hiptmair2018shape} and \cite[Thm. 3.2]{hettlich1995frechet}. Here we have replaced the Sommerfeld radiation condition with the exact absorbing boundary condition at $\partial B_R$ via the Dirichlet-\-to-\-Neumann map \cite[Sect. 2.6.3]{nedelec2001acoustic}.
    In \eqref{eq:shapeder_pde}, $ \Bn$ denotes the field with unit norm which is normal to the interface with radius $r$, whereas $\nabla_{\Gamma}$ and $\nabla_{\Gamma}\cdot$ are the surface gradient and divergence, respectively. We also remind that $D_R$ is as in Theorem \ref{thm:MoiolaSpence}. Since $\Be_{\rho}$ denotes the unit vector in radial direction, 
    $\delta r \Be_{\rho}\cdot \Bn$ is the variation of the interface with radius $r$ in normal direction. 


\smallskip
\underline{\textsl{Step 2}} We proceed to showing linearity and continuity of the map $\delta r\mapsto \delta u(r; \delta r)$, from $C^{2,1}_{per}(\Dphi)$ to $H^1(B_R\setminus \Gamma(\omega))$. Linearity is clear from \eqref{eq:shapeder_pde}, where $\delta r$ enters via the jump conditions. For continuity, it is then sufficient to show the boundedness of the map (see also proof of Proposition 3.5 in \cite{bui2014analysis}). Uniformity of the bound with respect to the radius will allow us to use the bound in a whole neighborhood of $r\in X$ in Step 3. We apply Theorem \ref{thm:MoiolaSpence2} to \eqref{eq:shapeder_pde} to bound
\begin{equation}\label{eq:deltau_bound}
     \begin{split}
     \lVert \delta u \rVert_{H^1_{\kappa_0,\alpha,n}(D_{in}\cup D_R)}^2\leq C_1 \lVert \nabla_{\Gamma} g_D\rVert^2_{L^2(\Gamma)} + C_2\lVert g_D\rVert^2_{L^2(\Gamma)} + C_3 \lVert g_N\rVert_{L^2(\Gamma)}^2, 
    \end{split}
\end{equation}
where the constants $C_1, C_2$ and $C_3$ can be bounded independently of $r\in X$ thanks to Assumptions \ref{ass:Expansion}-\ref{ass:uniform} and depend on the material coefficients and $\kappa_0$ as in \eqref{eq:normbound2}. Moreover,
\begin{equation*}
    g_D =-\delta r \Be_{\rho}\cdot \Bn\left\llbracket \frac{\partial u}{\partial \Bn}\right\rrbracket,  \quad   g_N=\nabla_{\Gamma} \cdot \left(\llbracket \alpha \rrbracket \delta r(\Be_{\rho}\cdot \Bn)\nabla_{\Gamma} (u+u^i) \right)+(\delta r\Be_{\rho}\cdot \Bn)\kappa_0^2\llbracket n \rrbracket (u+u^i).
\end{equation*}
Thanks to the regularity assumptions on the interface which ensure the solution $u$ to have $H^3$-regularity in each subdomain \cite[Ch. 4]{mclean2000strongly}, we can bound 
\begin{align*}
    \lVert g_D\lVert_{L^2(\Gamma)} &\leq \lVert \delta r \rVert_{C^0_{per}(\Dphi)}\lVert  \Bn\rVert_{C^0(\Gamma)}\Big\lVert\Big\llbracket \frac{\partial u}{\partial \Bn}\Big\rrbracket\Big\rVert_{L^2(\Gamma)},\\
    \lVert \nabla_{\Gamma} g_D\lVert_{L^2(\Gamma)} &\leq 3\lVert \delta r \rVert_{C^1_{per}(\Dphi)}\lVert  \Bn\rVert_{C^1(\Gamma)}\Big\lVert\Big\llbracket \frac{\partial u}{\partial \Bn}\Big\rrbracket\Big\rVert_{H^1(\Gamma)},\\
    \lVert g_N\lVert_{L^2(\Gamma)} &\leq \left(3\llbracket \alpha \rrbracket + \kappa_0^2 \llbracket n\rrbracket\right)\lVert \delta r \rVert_{C^1_{per}(\Dphi)}\lVert  \Bn\rVert_{C^1(\Gamma)}\lVert u + u^i\rVert_{H^2(\Gamma)},
\end{align*}
which, inserted in \eqref{eq:deltau_bound}, imply
\begin{equation}\label{eq:upperboundforlipschitz}
\begin{split}
    \lVert \delta u \rVert_{H^1_{\kappa_0,\alpha,n}(D_{in}\cup D_R)}^2\leq C \lVert \delta r \rVert_{C^1_{per}(\Dphi)}^2\lVert  \Bn\rVert_{C^1(\Gamma)}^2(\lVert u\rVert_{H^3(D_{in})}^2+\lVert u\rVert_{H^3(D_R)}^2 + \lVert u^i\rVert_{H^2(B_R)}^2), 
\end{split}
\end{equation}
for a constant $C=C(r)>0$ depending on the material coefficients, on $\kappa_0^2$ and on $r$ but not on its variation $\delta r$. Note that we can bound the $C^1$-norm of the radius variation with its $C^{2,1}$-norm.

In \eqref{eq:upperboundforlipschitz}, the quantities depending on $r\in X$ are the constant $C(r)$, whose dependence is due to the application of the trace inequality, and the piecewise $H^3$-norm of $u$. In view of the application of the mean value theorem in Step 3, we remark that these quantities can be bounded from above uniformly in $r\in X$. For the piecewise $H^3$-norm, a way to conclude this is to first consider a mapping approach to a reference geometric configuration (e.g., the one corresponding to the average radius $r_0$) and obtain bounds independent of $r\in X$ on the solution on the reference configuration using Theorems 4.16 and 4.20 in \cite{mclean2000strongly}, see for instance \cite{hiptmair2018large}. Then, the bound for the piecewise $H^3$-norm on the physical configuration follows from it by using the radius-independent upper bound on the piecewise $C^{2,1}$-norm of the mapping, in the same spirit as Lemma 1 in \cite{harbrecht2016analysis}. For $C(r)$, we need a uniform bound on the constant in the trace inequality, and, for this, explicit expressions as from \cite{auchmuty2014sharp} come to our help. Namely, from Theorem 6.1 in \cite{auchmuty2014sharp} we see that, if we can bound $\lVert u\rVert_{L^2(D_{in})}$, $\lVert u\rVert_{L^2(D_R)}$ and the $W^{1,\infty}$-norms of the solutions to a Poisson's equation with Neumann boundary conditions in $D_{in}$ and $D_R$ respectively (see Eq. (3.5) in \cite{auchmuty2014sharp}), then the constant in the trace inequality is bounded from above uniformly in $r\in X$. The bounds on the $L^2$-norms are a consequence of the bound on the piecewise $H^3$-norm, and our regularity assumptions on $\mu_0$-a.e. radius realization allow to bound the $W^{1,\infty}$-norms of the solutions to the Poisson's equations in both subdomains by using again a mapping approach and Schauder's regularity estimates \cite[Ch. 6]{gilbarg2015elliptic}.


\smallskip
\underline{\textsl{Step 3}} To apply the mean value theorem \cite[Thm. 3.2.7]{drabek2007methods}, we cannot work with a space for the solution $u$ which is dependent on the realization $r\in X$. This is why we only consider the restriction of $u$ to $B_R\setminus U$. We obtain that, for $\mu_0$-a.e. $r\in X$, there exists a ball $B_{\rho}(r)\subset C^{2,1}_{per}(\Dphi)$ centered in $r$ and with radius $\rho<\frac{1}{2}\inf_{\omega\in\Omega,\varphi_\in\Dphi}r(\omega,\varphi)$ such that, for all $r_1,r_2\in B_{\rho}(r)\cap X$,
\begin{equation*}
   \lVert u_1 - u_2 \rVert_{H^1_{\kappa_0,\alpha,n}(B_R\setminus U)} \leq L \lVert r_1-r_2\rVert_{C^{2,1}_{per}(\Dphi)},
\end{equation*}
where $u_1$ and $u_2$ are the solutions where the interfaces are determined by $r_1$ and $r_2$ respectively, $L$ is the upper bound on the right-hand side of \eqref{eq:upperboundforlipschitz}, and the two subdomains for the broken norm coinciding with the regions of $B_R$ which are either always inside or outside the scatterer. This shows the desired local Lipschitz continuity.

\end{proof}

The following result, for the case $\alpha_{in}=\alpha_{out}$, relaxes the smoothness assumptions on the scatterer's boundary and considers the range of the forward operator to be the whole $H^1(B_R)$.

\begin{proposition}\label{prop:helmholtz_continuity2}
    Assume that $\alpha_{out} = \alpha_{in}$ in Definition \ref{def:HSTP}, and that Assumptions \ref{ass:Expansion}-\ref{ass:uniform} hold. Furthermore, assume that 
    $\frac{n_{in}}{n_{out}}\leq 1$. Then the forward operator $G:X\rightarrow V$ mapping the radius to the solution to \eqref{eq:HSTP} on $B_R$, with $X$ as in \eqref{eq:radiusspace} and $V=H^1(B_R)$, is locally H\"older $\mu_0$-continuous.     
\end{proposition}
\begin{proof}
Let $u_1$ resp. $u_2$ denote the solution to equation \eqref{eq:HSTP} when the scatterer boundary is parametrized by the radius $r_1$ resp. $r_2$. We introduce the notation $\alpha_1, n_1$ for the piecewise constant coefficients \eqref{eq:alphakappacoeffs} jumping at the interface parametrized by $r_1$, and similarly for $\alpha_2,n_2$. Moreover, we denote $D^{1\backslash2}\coloneqq D_{in,1}\backslash D_{in,2}$ and $D^{2\backslash1}\coloneqq D_{in,2}\backslash D_{in,1}$, where $D_{in,1}$ resp. $D_{in,2}$ are the part of the domain inside the scatter in the case of the interface being parameterized by $r_1$ resp. $r_2$ (and similarly later for $D_{out,1}, D_{out,2}$). Observe that $u_1-u_2$ is already a scattering field, namely it satisfies \eqref{eq:HSTP2_sommerfeld}. 
Since $\alpha_{out}=\alpha_{in}$, we then find that $(u_1-u_2)\in H^1(B_R)$ satisfies the Helmholtz transmission problem in Theorem \ref{thm:MoiolaSpence} with $\alpha_1$ and $n_1$ on the right-hand side, $f_{in}=\kappa_0^2(n_{out}-n_{in})(u_2+u^i)\mathbbm{1}_{D^{1\backslash2}}$, $f_{out} =\kappa_0^2(n_{in}-n_{out})(u_2+u^i)\mathbbm{1}_{D^{2\backslash1}}$, where $\mathbbm{1}_D$ denotes the indicator function on a domain $D$, and $g_N=g_D\equiv 0$. The bound \eqref{eq:normbound}  applied to the difference $(u_1-u_2)$ gives us:
\begin{equation*}
    \begin{split}
        & \lVert u_1-u_2 \rVert_{H^1_{\kappa_0,\alpha,n}(D_{in,1}\cup D_{out,1})}^2 \\
        &\leq \left[\frac{4 \text{diam}(D_{in,1})^2}{\alpha_{in}}+\frac{1}{n_{in}}\left(2\sqrt{\frac{n_{out}}{\alpha_{out}}}R+\frac{1}{\kappa_0} \right)^2\right]\left(\kappa_0^4(n_{out}-n_{in})^2\|u_2+u^i\|^2_{L^2(D^{1\backslash2})}\right)\\
        & + \left[\frac{4 R^2}{\alpha_{out}}+\frac{1}{n_{out}}\left(2\sqrt{\frac{n_{out}}{\alpha_{out}}}R+\frac{1}{\kappa_0} \right)^2\right]\left(\kappa_0^4(n_{in}-n_{out})^2\|u_2+u^i\|^2_{L^2(D^{2\backslash1})}\right).
    \end{split}
\end{equation*}
We can use H\"older's inequality to get $\|(u_2+u^i)\|^2_{L^2(D^{1\backslash2})}\leq \sqrt{\left|D^{1\backslash2}\right|}\|(u_2+u^i)^2\|_{L^2(D^{1\backslash2})}$, 
and, since $u_2, u^i\in H^1(B_R)$, the Sobolev embedding theorem \cite[Sect. 7.7]{gilbarg2015elliptic} ensures that the right-hand side is finite for $d=2,3$. An analogous result can be obtained for the $L^2$-norm on $D^{2\backslash1}$. The claim follows by noticing that $\left|D^{1\backslash2}\right|$ and $\left|D^{2\backslash1}\right|$ are bounded from above by $\lVert r_1-r_2\rVert_{C^0_{per}(\Dphi)}$ (with a constant dependent on the dimension $d$ and on $\max\left\{\lVert r_1\rVert_{C^0_{per}(\Dphi)},\lVert r_2\rVert_{C^0_{per}(\Dphi)}\right\}$).
\end{proof}

\begin{remark}[Lipschitz continuity for globally constant $\alpha$]
    When $\alpha_{in}=\alpha_{out}$ and the radius is only $\mu_0$-a.s. Lipschitz continuous, the last proposition allows us to conclude the measurability of the forward operator via H\"older $\mu_0$-continuity, which is sufficient for the well-posedness of our Bayesian inverse problem. The more general setting of Proposition \ref{prop:helmholtz_continuity1} requires instead $C^{2,1}$-smoothness for the radius and a slightly more restrictive space $V$, but it allows to conclude Lipschitz rather than only H\"older continuity. 
    A close inspection of the proof of Proposition \ref{prop:helmholtz_continuity1} reveals although that, when $\alpha_{in}=\alpha_{out}$, it is sufficient to ask for a $\mu_0$-a.s. $C^{1,1}$-continuous radius to conclude Lipschitz continuity of the forward operator from $X$ to $V=H^1(B_R)$: when $\alpha_{in}=\alpha_{out}$, the jump terms $\llbracket \alpha\rrbracket$ and $\llbracket \frac{\partial u}{\partial \Bn}\rrbracket$ vanish (cf. \eqref{eq:jump2}), and those are the terms that require stricter smoothness requirements and hinder global $H^1$-smoothness for the shape derivative.
\end{remark}

We now address a wavenumber-explicit bound on the $L^{\infty}$-norm of the forward operator, which, according to $(ii)$ in Corollary \ref{cor:pdewellposedness}, will give us a wavenumber-explicit stability bound for the posterior dependence on the data. We still denote by $B_R$ a ball with a radius $R$ large enough to contain all scatterer's realizations. The following is then a trivial consequence of Corollary \ref{cor:uscattbound}.

\begin{proposition}\label{prop:helmholtz_stability}
  Let Assumptions \ref{ass:Expansion}-\ref{ass:uniform} hold, and let the material parameters in \eqref{eq:alphakappacoeffs} satisfy \eqref{eq:nontrappingcond}.   
    Then, for the forward operator $G:X\rightarrow V$ mapping the radius to the solution to \eqref{eq:HSTP}, with $X$ as in \eqref{eq:radiusspace} and either $V=H^1(B_R)$ or $V=H^1(B_R\setminus U)$ equipped with the norm \eqref{eq:weightednorm}, we have
    \begin{equation}
    \begin{split}
                &\lVert G\rVert_{L^{\infty}_{\mu_0}(X,V)}\\
                &\leq C_{\kappa_0} C_{1}\alpha_{out}\lVert \nabla u^i\rVert_{L^2(D_{out,H})} + (C_{\kappa_0} C_{2}\alpha_{out} + \sqrt{\alpha_{out}} C_{1})\lVert u^i\rVert_{L^2(D_{out,H})} + \lVert u^i \rVert_{H^1_{\kappa_0,\alpha,n}},
    \end{split}
    \end{equation}
    where $C_{\kappa_0}$, $C_1$ and $C_2$ are as in Corollary \ref{cor:uscattbound}. The domain $D_{out,H}$ is as in \eqref{eq:holdall}, and in $C_1,C_2$ we take $R_{scatt}>(1+\gamma_{\beta})\lVert r_0\rVert_{C^0_{per}(\Dphi)}$ with $\gamma_{\beta}$ as in Assumption \ref{ass:uniform}.
\end{proposition}

Summarizing Propositions \ref{prop:helmholtz_continuity1}-\ref{prop:helmholtz_stability} in the framework of Corollary \ref{cor:pdewellposedness}, we get our main result.
\begin{theorem}\label{thm:wellposed_HSTP}
    Let Assumptions \ref{ass:Expansion} and \ref{ass:uniform} hold, and assume one of the following settings:
    \begin{itemize}
        \item 
        $\alpha_{in}=\alpha_{out}$, $\frac{n_{in}}{n_{out}}<1$ and $V=H^1(B_R)$;
        \item the radius $r$ in \eqref{eq:radiusomega} is $\mu_0$-a.s. of class $C^{2,1}$ and $\lVert r\rVert_{C^{2,1}_{per}(\Dphi)}$ is uniformly bounded on $X$,
         $\frac{n_{in}}{n_{out}}<1<\frac{\alpha_{in}}{\alpha_{out}}$ and $V=H^1(B_R\setminus U)$, with $U$ as in \eqref{eq:u}.
    \end{itemize}
Furthermore, assume, in either of the two cases, that the observation operator is continuous on $V$, namely $\mathcal{O}=(\mathcal{O}_k)_{k=1}^K\in (V^{\ast})^K$, and that the noise is Gaussian, such that the potential is given by \eqref{eq:potential} for a symmetric positive definite matrix $\Sigma$ with minimum eigenvalue $\lambda_{min}$. Then the Bayesian shape inverse problem stated in Subsection \ref{ssec:setting} is well-posed, namely:
   \begin{enumerate}
       \item[(i)] the posterior measure is absolutely continuous with respect to the prior and has Radon-Nikodym derivative \eqref{eq:RadonNikodym};
       \item[(ii)] the posterior depends continuously on the data: for each $\gamma>0$ such that $|\boldsymbol{\delta}|$, $|\boldsymbol{\delta}'|\leq \gamma$,
       \begin{equation}\label{eq:dHell_Helmholtz}
       	\begin{split}
       		d_{\mathrm{Hell}}(\mu^{\boldsymbol{\delta}},\mu^{\boldsymbol{\delta'}})&\leq \Tilde{C}_{\mathrm{Hell}}C_{\gamma,G}\exp\left(\frac{1}{2}C_{\gamma,G}^2\right)|\boldsymbol{\delta} - \boldsymbol{\delta}'|,\\
       		d_{W^1}(\mu^{\boldsymbol{\delta}},\mu^{\boldsymbol{\delta'}})&\leq\Tilde{C}_{W^1} C_{\mu_0}C_{\gamma,G}\exp\left(C_{\gamma,G}^2\right)|\boldsymbol{\delta} - \boldsymbol{\delta}'|,
       	\end{split}
       \end{equation}
       for some $\Tilde{C}_{\mathrm{Hell}}, \Tilde{C}_{W^1}>0$ independent of $\gamma$, the wavenumber, and any other physical parameter entering the forward model, $C_{\mu_0}=\sum_{j\geq 1} |\beta_j|$ and      
       \begin{equation}\label{eq:stabconst_helmholtz}
       	\begin{split}
       		C_{\gamma,G} &= \lambda_{min}^{-1}\gamma
       		+ \lambda_{min}^{-1}\Onorm \left[ \lVert u^i \rVert_{H^1_{\kappa_0,\alpha,n}}+ C_{\kappa_0} C_{1}\alpha_{out}\lVert \nabla u^i\rVert_{L^2(D_{out,H})} \right.\\
       		& + \left.\left(C_{\kappa_0} C_{2}\alpha_{out} + \sqrt{\alpha_{out}} C_{1}\right)\lVert u^i\rVert_{L^2(D_{out,H})} \right].
       \end{split}
       \end{equation}
       Here $C_{\kappa_0}$, 
    $C_1$ and $C_2$ are as in Corollary \ref{cor:uscattbound} with $R_{scatt}>(1+\gamma_{\beta})\lVert r_0\rVert_{C^0_{per}(\Dphi)}$ and $\gamma_{\beta}$ as in Assumption \ref{ass:uniform}.
   \end{enumerate}
\end{theorem}
To help the intuition in the stability estimate \eqref{eq:stabconst_helmholtz}, we consider the case of $u^i$ being a plane wave. One can then see that the behavior with respect to the wavenumber is characterized by
       \begin{equation}\label{eq:stabconst_helmholtz_pw}
         C_{\gamma,G}\sim \lambda_{min}^{-1}\gamma + \lambda_{min}^{-1}\left|\left(\lVert\mathcal{O}_k\rVert_{V^{\ast}}\right)_{k=1}^K\right| \kappa_0 R,
       \end{equation}
where $\sim$ hides constants depending on the material and geometrical parameters but not $\kappa_0$.
Using this and the theorem's statement, we can then quantify facts, some of which might be expected from a physical perspective:
\begin{itemize}
    \item As it is well-known from previous results on the forward model \cite{moiola2019acoustic, hiptmair2025frequency}, the lengthscale of the problem is dictated by the wavenumber $\kappa_0$, meaning that what matters most for stability is the size of the domain relative to the wavenumber, expressed by the term $\kappa_0 R$.
    \item The geometric quantities $R_{scatt}$ and $R$ do not only depend on the physical setting of the scattering problem, but also on the inverse problem assumptions and setting. Namely, $R_{scatt}$ also depends on how much shape variation we allow in the prior, and $R$ on the experimental setup, in the sense that $R$ should be large enough such that the observation operator is continuous; e.g., for observations consisting of smoothed point values of the field, the farther the measurements, the larger $R$ needs to be. The larger $R-R_{scatt}$, which could correspond to measurements taken farther away, the larger the stability constant. 
    \item When $\alpha_{in}=\alpha_{out}$, not only we can relax the smoothness assumptions on the scatterer's boundary, but we can allow for a larger set of observation operators. For instance, it allows for localized measurements very close to the scatterer.
    \item The stability estimate \eqref{eq:dHell_Helmholtz} tells us how much the posterior is sensitive to noise, coming not only from measurement error but also possibly from some modeling approximation (the so-called model discrepancy).  Indeed, considering \eqref{eq:dHell_Helmholtz} when $\Bdelta=\mathcal{G}(r^{\dag})$ and $\Bdelta'=\mathcal{G}(r^{\dag})+\eta$, with $r^{\dag}$ an underlying \textsl{truth} and $\left|\eta\right|=\gamma$, the effect of the noise is expressed by $\left|\Bdelta-\Bdelta'\right|=\gamma$. Our estimates tell us that this is amplified by a factor $\sim\kappa_0 R$: for a fixed geometrical configuration (and with the same signal to noise ratio, see Subsection \ref{ssec:numres}), the higher the frequency, the larger the sensitivity. 
\end{itemize}

\begin{remark}[Role of material contrast]
    The estimate \eqref{eq:dHell_Helmholtz} is explicit and optimal in the wavenumber, but less explicit in the material parameters \cite[Rmk. 3.6]{chaumont2023explicit}. The effect of the contrast between the two materials on stability is clearer from the (wavenumber-\-suboptimal) bound on the scattered wave as from Corollary 3.1 in \cite{moiola2019acoustic}, which follows from Theorem \ref{thm:MoiolaSpence2} applied to the scattered wave. The resulting constant for \eqref{eq:stabconst_helmholtz} would then be
    \begin{align*}
         C_{\gamma,G}&=\lambda_{min}^{-1}\gamma + \lambda_{min}^{-1}\left|\left(\lVert\mathcal{O}_k\rVert_{V^{\ast}}\right)_{k=1}^K\right|\left[\frac{4(\kappa_0 \diam(D_{in,H}))^2}{\alpha_{in}} + \frac{(\kappa_0 R)^2}{n_{in}}\left(2\sqrt{\frac{n_{out}}{\alpha_{out}}} + \frac{d-1}{\kappa_0 R}\right)^2\right]^{\frac{1}{2}}\\
         &\cdot \kappa_0\left(\frac{\alpha_{in}}{\alpha_{out}}n_{in}-n_{out}\right)\lVert u^i\rVert_{L^2(D_{in,H})}\\
     &+ \lambda_{min}^{-1}\left|\left(\lVert\mathcal{O}_k\rVert_{V^{\ast}}\right)_{k=1}^K\right|\left[\frac{\diam(D_{in,H})(4\alpha_{in}+2\alpha_{out})}{\alpha_{in}-\alpha_{out}}+\frac{2\left(n_{out}(\kappa_0 R)^2+\alpha_{out}\frac{(d-1)^2}{4}\right)}{\kappa_0^{2}\diam(D_{in,H}\setminus U))(n_{out}-n_{in})}\right]^{\frac{1}{2}}\\
     &\cdot \frac{2 C_{surf}(\alpha_{in}-\alpha_{out})}{\hat{\gamma}\alpha_{out}}\lVert u^i\rVert_{C^1(U)}, 
    \end{align*}
where $U$ is as in \eqref{eq:u}, $C_{surf}:=\sup_{\omega\in\Omega}|\Gamma (\omega)|^{\frac{1}{2}}$ is bounded from above by a constant depending on $\operatorname{ess}\sup_{r\in X}\lVert r\Vert_{C^{0,1}_{per}(\Dphi)}$ only, and $\hat{\gamma}$ is the minimum between $\frac{1}{2}$ and $\tilde{\gamma}$ in Lemma \ref{lem:starshaped}.
\end{remark}



\subsection{Extension to Sound-Soft Obstacle Scattering}\label{ssec:helmholtzextensions}

Thanks to the estimates in \cite{graham2019helmholtz}, which extend those in \cite{moiola2019acoustic}, we can also extend our results to the exterior Dirichlet problem with heterogeneous low-order coefficient. Referring to Figure \ref{fig:particle}, this corresponds to a not penetrable scatterer acting as an obstacle, such that 
the total field vanishes at its boundary. This is called sound-soft scattering in acoustics \cite[Sect. 2.1]{colton1998inverse}. As in the transmission problem, $D_{out}$ is the exterior of the scatterer, $\Gamma$ its boundary parametrized by \eqref{eq:radiusomega} with normal $\Bn$ and we use the notation \eqref{eq:holdall} for the hold-all domains.

\begin{definition}\label{def:HSSP} (Sound-soft scattering problem in heterogeneous medium)\\
    Let $\kappa_0 \in \mathbb{R}_{>0}$ be the wavenumber, and let $n=n(\Bx)$ be a real-valued function defined on $D_{out,H}$, bounded away from zero and such that $1-n$ has compact support. Let $u^i$ be a solution of $\Delta u^i+\kappa_0^2u^i=0$, that is $C^\infty$ in a neighbourhood of $\overline{D_{in,H}}$. For $\omega\in\Omega$, a solution $u=u(\omega)$ to the sound-soft scattering problem solves
    \begin{subequations}\label{eq:HSSP}
\begin{empheq}[left={\empheqlbrace}]{alignat=2}
        &-\Delta u - \kappa_0^2n(\Bx)u =f\quad && \text{in } D_{out}(\omega),\\
        & u = u^i \quad && \text{on }\Gamma(\omega),\\   
        &\lim_{\lVert\Bx\rVert\rightarrow\infty}\lVert\Bx\rVert^{\frac{d-1}{2}}\left(\frac{\partial}{\partial \lVert\Bx\rVert}-i\kappa_0\right)u = 0,&&
\end{empheq}
\end{subequations}
with $f=\kappa_0^2 (1-n)u^i$.
\end{definition}
From \cite{graham2019helmholtz}, we will use mainly Theorem 2.19 $(ii)$. We need to assume some regularity, boundedness, and a non-trapping condition on the refraction index $n$. We refer to \cite[Sect. 7]{graham2019helmholtz} for a detailed physical interpretation of the latter.

\begin{assumption}[Condition 2.18 in \cite{graham2019helmholtz}]\label{ass:nDirichlet} The refraction index satisfies $n\in C^{0,1}(\overline{D_{out,H}})$, and there exist $n_{min},n_{max}>0$ such that $n_{min}\leq n(\Bx)\leq n_{max}<\infty$ for a.e. $\Bx\in D_{out,H}$. Furthermore, there exists $\mu_n>0$ such that $2n(\Bx)+\Bx\cdot\nabla n(\Bx)\geq \mu_n$ for a.e. $\Bx\in D_{out,H}$.
\end{assumption}

Introducing the set
\begin{equation}\label{eq:ud}
 U_D:=\left\{\Bx\in B_R: \left|\Bx\right|\leq r^{+}(\varphi_{\Bx})\right\},
\end{equation}
where $B_R$ 
contains all scatterer's realizations, we  
obtain the analogous of Proposition \ref{prop:helmholtz_continuity1}.
\begin{proposition}\label{prop:soundsoft_continuity}
    Under Assumptions \ref{ass:Expansion}-\ref{ass:uniform}, assume that the radius $r$ in \eqref{eq:radiusomega} is $\mu_0$-a.s. of class $C^{2,1}$ and that its $C^{2,1}$-norm is uniformly bounded on $X$, i.e. that $\mu_0(X\cap C^{2,1}_{per}(\Dphi))=1$ and $\operatorname{ess} \sup_{r\in X}\lVert r \rVert_{C^{2,1}_{per}(\Dphi)}<\infty$.
    Furthermore, assume that the refraction index fulfills Assumption \ref{ass:nDirichlet}, that the support of $1-n$ is contained in $B_R$ and $R\geq \sqrt{3/8}\kappa_0^{-1}$.
    
    Then the forward operator $G:X\rightarrow V$ mapping the radius to the solution to \eqref{eq:HSSP} on $B_R\setminus U_D$, with $X$ as in \eqref{eq:radiusspace} and $V=H^1(B_R\setminus U_D)$, is locally Lipschitz $\mu_0$-continuous. 
\end{proposition}
\begin{proof}
Here we only mention the main differences to the proof of Proposition \ref{prop:helmholtz_continuity1} in Steps 1 and 2. 

\medskip
\underline{\textsl{Step 1}} The $\mu_0$-a.s. $H^2$-regularity of the solution $u$ to \eqref{eq:HSSP} can also be obtained using the regularity results in \cite[Ch.4]{mclean2000strongly}. 
Corollary 3.1 in \cite{hiptmair2018shape} should now be applied in conjunction with Theorem 4.2 in \cite{hiptmair2018shape}. It follows that, for $r\in X$ and $\delta r\in C^{2,1}_{per}(\Dphi)$ with sufficiently small norm, the shape derivative $\delta u(r;\delta r)$ is the solution in $H^1(D_R(\omega))$, $D_R(\omega)=D_{out}(\omega)\cap B_R$, of 
    \begin{equation}\label{eq:shapeder_pde_ss}
        \begin{cases}
        -\Delta \delta u - \kappa_0^2n(\Bx)\delta u=0\quad&\text{in }D_R(\omega),\\
         \delta u = -\delta r\Be_{\rho}\cdot \Bn\dfrac{\partial (u+u^i)}{\partial \Bn} \quad &\text{on }\Gamma(\omega),\\   
         \dfrac{\partial \delta u}{\partial \Bn}=\operatorname{DtN}(\delta u) \quad &\text{on }\partial B_R,
    \end{cases}
    \end{equation}
    with boundary condition to be understood in a weak sense \cite{hiptmair2018shape} and $u$ solving \eqref{eq:HSSP}, see also \cite{bui2014analysis}.

\medskip
\underline{\textsl{Step 2}} Theorem \ref{thm:MoiolaSpence2} should be replaced by Theorem 2.19 $(ii)$ in \cite{graham2019helmholtz}, yielding
\begin{equation}
    \mu_n\left(\lVert\nabla \delta u\rVert^2_{L^2(B_R)}+\kappa_0^2 \lVert\delta u\rVert^2_{L^2(B_R)}\right)\leq C_1 \lVert \nabla_{\Gamma} g_D \rVert^2_{L^2(\Gamma)} + C_2\kappa_0^2 \lVert g_D \rVert^2_{L^2(\Gamma)},
\end{equation}
where $g_D = -\delta r \Be_{\rho}\cdot \Bn\dfrac{\partial (u+u^i)}{\partial \Bn}$ and $C_1,C_2$ can be bounded from above independently of $r\in X$ (by taking the maximum possible diameter of the scatterer). The terms $\lVert \nabla_{\Gamma} g_D \rVert^2_{L^2(\Gamma)}$ and $\lVert g_D \rVert^2_{L^2(\Gamma)}$ can then be bounded analogously to the proof of Proposition \ref{prop:helmholtz_continuity1}, with Theorem 4.20 in \cite{mclean2000strongly} to be replaced by Theorem 4.18 in \cite{mclean2000strongly}.
%
\end{proof}

The next result, 
analogue to Proposition \ref{prop:helmholtz_stability}, follows by applying the bound from Theorem 2.19 $(ii)$ in \cite{graham2019helmholtz} to \eqref{eq:HSSP}. In the statement below, the weighted norm 
 is defined as
\begin{equation}\label{eq:weightednorm2}
        \lVert v \rVert_{H^1_{\kappa_0}(B_R\setminus U_D)}^2:=  \lVert \nabla v\rVert_{L^2(B_R\setminus U_D)}^2 + \kappa_0^2\lVert v\rVert^2_{L^2(B_R\setminus U_D)}.
\end{equation}

\begin{proposition}\label{prop:soundsoft_stability}
  Let Assumptions \ref{ass:Expansion}-\ref{ass:uniform} and \ref{ass:nDirichlet} hold, and let the support of $1-n$ be contained in $B_R$ with $R\geq \sqrt{3/8}\kappa_0^{-1}$. Then, for the forward operator $G:X\rightarrow V$ mapping the radius to the solution to \eqref{eq:HSSP} in $B_R\setminus U_D$, with $X$ as in \eqref{eq:radiusspace} and $V=H^1(B_R\setminus U_D)$ equipped with the norm \eqref{eq:weightednorm2}, we have
    \begin{equation}\label{eq:boundedness_ss}
        \lVert G\rVert_{L^{\infty}_{\mu_0}(X,V)}\leq C_1 \kappa_0 \lVert (1-n)u^i\rVert_{L^2(D_{out,H})} + C_2 \lVert u^i\rVert_{C^1(U_D)} + C_3\lVert u^i\rVert_{C^0(U_D)},
        \end{equation}
        with $D_{out,H}$ as in \eqref{eq:holdall}, $U_D$ as in \eqref{eq:ud} and
\begin{equation}\label{eq:constants_ss}
\begin{split}
    C_1 &= 2\left(\frac{4(\kappa_0 R)^2 }{\mu_n^2}\left(1 + \left(2+\frac{d-2}{2\kappa_0 R}\right)^2\right)\left(1+\frac{3}{2}n_{max}\right)^2 + \frac{2}{n_{max}}\right)^{\frac{1}{2}},\\
    C_2 &= C_{surf} \sqrt{\frac{2}{\mu_n}}\left(1+\frac{3}{2}n_{max}\right)^{\frac{1}{2}}(\diam (D_{in, H}))^{\frac{1}{2}}\left(1+\frac{4 \diam (D_{in, H})}{\tilde{\gamma}}\right)^{\frac{1}{2}},\\
    C_3 &= 2 C_{surf}\left(\frac{8}{\mu_n}\left(1+\frac{3}{2}n_{max}\right)\frac{(\kappa_0 R)^2}{\tilde{\gamma}}\left(2+\frac{d-2}{2\kappa_0 R}\right)^2 +\frac{2}{\tilde{\gamma}}\right)^{\frac{1}{2}}.
\end{split}
    \end{equation}
Here $C_{surf}=\sup_{\omega\in\Omega}|\Gamma (\omega)|^{\frac{1}{2}}$ is bounded from above by a constant only depending on $\operatorname{ess}\sup_{r\in X}\lVert r\Vert_{C^{0,1}_{per}(\Dphi)}$, and $\tilde{\gamma}$ is as in Lemma \ref{lem:starshaped}. Moreover, for $\gamma_{\beta}$ as in Assumption \ref{ass:uniform}, $\diam(D_{in,H})\!\leq\! 2(1+\gamma_{\beta})\lVert r_0\rVert_{C^0_{per}(\Dphi)}$.
\end{proposition}
In the estimates above, the amount of shape variation in the prior affects the domain of integration in the norms. Although $C_2$ contains quantities related to the size of the domain in absolute values, it is multiplied with the $C^1$-norm of the incoming wave, which often contains a multiplication by $\kappa_0$ (for instance, if $u^i$ is a plane wave), leading again to a relative importance of size with respect to the wavenumber. Similarly, the term multiplied by $C_1$ involves integration over a subdomain of $D_{out,H}$, due to $(1-n)$ having compact support, and the $L^2$-norm multiplied by $\kappa_0$ entails an interplay between the wavenumber and 
the size of the support of $(1-n)$.

In \eqref{eq:boundedness_ss}-\eqref{eq:constants_ss}, the realization-independent upper bounds for the second and third summands have been obtained from Theorem 2.19 $(ii)$ in \cite{graham2019helmholtz} by bounding trace norms 
with the $C^0$- and $C^1$-norms of $u^i$ on the volume. 
We did not use trace inequalities to avoid higher smoothness requirements on the radius and because the resulting bounds would involve an $H^2$-norm of $u^i$, which, 
 for instance when $u^i$ is a plane wave, would be suboptimal 
 in terms of $\kappa_0$-dependence.

The previous results with Corollary \ref{cor:pdewellposedness} yield well-posedness for sound-soft scattering.

\begin{theorem}\label{thm:wellposed_soft}
    Under Assumptions \ref{ass:Expansion}-\ref{ass:uniform}, assume that the radius $r$ in \eqref{eq:radiusomega} is $\mu_0$-a.s. of class $C^{2,1}$ and $\lVert r\rVert_{C^{2,1}_{per}(\Dphi)}$ is uniformly bounded on $X$, that the refraction index fulfills Assumption \ref{ass:nDirichlet}, and that the support of $1-n$ is contained in $B_R$ with $R\geq \sqrt{3/8}\kappa_0^{-1}$.
Furthermore, assume that the observation operator is continuous on $V=H^1(B_R\setminus U_D)$, namely $\mathcal{O}=(\mathcal{O}_k)_{k=1}^K\in (V^{\ast})^K$, and that the noise is Gaussian, such that the potential is given by \eqref{eq:potential} for a symmetric positive definite matrix $\Sigma$ with minimum eigenvalue $\lambda_{min}$. Then the Bayesian shape inverse problem stated in Subsection \ref{ssec:setting} when $u$ in \eqref{eq:AdditiveNoise} solves \eqref{eq:HSSP} is well-posed, namely:
   \begin{enumerate}
       \item[(i)] the posterior measure is absolutely continuous with respect to the prior and has Radon-Nikodym derivative \eqref{eq:RadonNikodym};
       \item[(ii)] the posterior depends continuously on the data: for each $\gamma>0$ such that $|\boldsymbol{\delta}|$, $|\boldsymbol{\delta}'|\leq \gamma$, we have \eqref{eq:dHell_Helmholtz} where now, for $C_1,C_2,C_3$ as in \eqref{eq:constants_ss},
       \end{enumerate}
       \begin{equation}\label{eq:stabconst_helmholtzsoft}
       \begin{split}
         C_{\gamma,G}&=\lambda_{min}^{-1}\gamma  \\
         & + \lambda_{min}^{-1}\left|\left(\lVert\mathcal{O}_k\rVert\right)_{k=1}^K\right|\left(C_1 \kappa_0 \lVert (1-n)u^i\rVert_{L^2(D_{out,H})} + C_2 \lVert u^i\rVert_{C^1(U_D)} + C_3\lVert u^i\rVert_{C^0(U_D)}\right).
       \end{split}
       \end{equation}
\end{theorem}

This result allows for similar considerations as those we made for the transmission problem.

\section{Numerical Illustrations}\label{sec:Results}
We provide, via some numerical experiments, further qualitative insight on the analytical results. 

The numerical methods are described in Subsection \ref{ssec:finiteelements} and Subsection \ref{ssec:smc}, on the forward solver and the sampling algorithm, respectively. Subsection \ref{ssec:settings} specifies the values for algorithmic and physical parameters used in our simulations and Subsection \ref{ssec:numres} presents numerical results for the inversion using two different frequencies, to see how the wavenumber affects the posterior sensitivity to noise. We focus on the transmission problem \eqref{eq:HSTP} in two dimensions ($d=2$).

All algorithms have been implemented in Python, using DOLFINx to implement the forward solver \cite{alnaes2014unified, baratta2023dolfinx,scroggs2022construction,scroggs2022basix} and starting from the code offered in \cite{bulte2020practical} to implement the sampling.


\subsection{Numerical Solution of the Forward Model}\label{ssec:finiteelements}
%
We base our discretization of \eqref{eq:HSSP} on the finite element method. For obstacle scattering, when the measurements occur at locations which are outside every scatterer realization, the boundary element method as considered in \cite{dolz2023solving} is a valid alternative.

We consider the truncated domain $B_R$ as in Section \ref{sec:well-posedness} and at $\partial B_R$ we approximate the Sommerfeld radiation condition \eqref{eq:Sommerfeld} with impedance boundary conditions: 
\[
\dfrac{\partial u}{\partial \Bn_R} - i\kappa_0 \sqrt{\frac{n_{out}}{\alpha_{out}}} u = 0\quad\text{on }\partial B_R,
\]
where $\Bn_R$ denotes the outer normal to $\partial B_R$.

Rather than discretizing \eqref{eq:HSSP} directly, we consider a mapping to a reference configuration \cite{xiu2006numerical,tartakovsky2006stochastic}, corresponding to the geometrical configuration where the boundary of the scatterer is described by the radius $r_0$. As domain mapping, we use radial rescaling with a piecewise linear mollifier whose support is contained inside $B_R$, such that it is the identity in a neighborhood of $\partial B_R$, see \cite[pp. 35-36]{kuijpers2023bayesian} for its analytical expression. 

The resulting variational formulation on the reference domain, see \cite{hiptmair2018large} and \cite[p. 37]{kuijpers2023bayesian}, is discretized using globally continuous, piecewise polynomial finite elements on an unstructured, simplicial mesh. We will adapt the polynomial order to the wavenumber, as further discussed in Subsection \ref{ssec:numres}. Since the mollifier for the domain mapping is only piecewise smooth, we resolve with the mesh the curve along which it has a kink in order to preserve the full finite element convergence rate. We solve for the total wave and subtract the incident field to obtain measurements of the scattered wave. For the latter, we consider smoothed point evaluations: 
\begin{equation}\label{eq:smoothedpointval}
	\mathcal{O}_k(u) = \frac{1}{\sqrt{2\pi}\sigma_o}\int_{B_R} \exp\left(-\frac{\lVert \Bx-\Bx_k\rVert^2}{2\sigma_o^2}\right) u(\Bx)\dd\Bx,\quad k=1,\ldots,K,
\end{equation}
where $\sigma_o>0$ is the level of smoothing. We will choose the measurement locations $\left\{\Bx_k\right\}_{k=1}^K$ to be sufficiently far from the scatterer such that the domain mapping is the identity around them. When approximating \eqref{eq:smoothedpointval} by quadrature, for computational efficiency we discard the quadrature points where the Gaussian function is essentially zero (in our code, when smaller than $10^{-8}$).

\subsection{Sampling Algorithm}\label{ssec:smc}
To sample from the posterior, we use the sequential Monte Carlo (SMC) method. 
Since, for the uncertain radius, we will truncate the expansion after $J\in\mathbb{N}$ terms and we use the pushforward of the distribution of $\BY:=\left(Y_j\right)_{j=1}^J$ under \eqref{eq:radiusomega}, here we explain the algorithm in terms of samples for $\BY$. Accordingly, $\mu_{\BY,0}=$ $\mathcal{U}([-1,1]^J)$ denotes the prior on $\BY$ and $\mu^{\Bdelta}_{\BY}$ its posterior. The latter is as in \eqref{eq:posterior} when replacing $\mu_0$ by $\mu_{\BY,0}$ (and composing $\mathcal{G}$ with the map given by \eqref{eq:radiusomega}). We focus on the main aspects of the algorithm and we refer to \cite{doucet2001introduction,chopin2002sequential, beskos2015sequential} for a more comprehensive treatment. 

The general idea is to approximate a distribution with a swarm of weighted particles (i.e., samples), and to transition from the prior to the posterior via a sequence of intermediate distributions. We define these intermediate distributions via tempering \cite{beskos2015sequential}: for a given sequence of temperatures $0=T_0<T_1<\ldots<T_L=1$, we set
\begin{equation}\label{eq:mul}
 \mu^{\Bdelta}_{\BY,l}(d \By) \propto \exp{\left(-T_l \frac{1}{2}|\boldsymbol{\delta}-\mathcal{G}(r(\By))|^2_{\Sigma}\right)}\mu_{\BY,0}(d \By), \quad l=0,\ldots,L,
\end{equation}
such that $\mu^{\Bdelta}_{\BY,0}$ is the prior and $\mu^{\Bdelta}_{\BY,L}$ the posterior. 
Using \eqref{eq:mul}, we have the update 
\begin{equation*}
 \mu^{\Bdelta}_{\BY,l+1}(\dd \By) \propto \exp{\left(-(T_{l+1}-T_l) \frac{1}{2}|\boldsymbol{\delta}-\mathcal{G}(r(\By))|^2_{\Sigma}\right)}\mu_{\BY,l}(\dd \By), \quad l=0,\ldots,L-1.
\end{equation*}

In the algorithm, we first initialize $M$ particles $\left\{\BY_i^0\right\}_{i=1}^M$ by sampling from the prior and we set their weights to $W^0_i=\frac{1}{M}$, for $i=1,\ldots,M$. We denote by $\BY_i^l$ and $W_i^l$ the position and weight of the $i$-th particle at the $l$-th step, respectively. For $l=0,\ldots,L-1$, we then move from a particle approximation of $\mu^{\Bdelta}_{\BY,l}$ to $\mu^{\Bdelta}_{\BY,l+1}$ via:
\begin{enumerate}
    \item \textsl{Update step}: we update the distribution by updating and renormalizing the weights only:
     $$w_i^{l+1}=\exp{\left(-(T_{l+1}-T_l) \frac{1}{2}|\boldsymbol{\delta}-\mathcal{G}(r(\BY_i^l))|^2_{\Sigma}\right)}w_i^{l},\quad W_i^{l+1}=\frac{w_i^{l+1}}{\sum_{l=1}^M w_i^{l+1}}, \quad \text{for }i=1,\ldots,M.$$ 
Note that this step does not require any new evaluation of the forward model.
    \item \textsl{Resampling Step}: in order to improve the approximation in the high-probability regions, we resample the particle population with replacement, according to the weights $W_i^{l+1}$, $i=1,\ldots,M$. After resampling, the weights are updated to $W_i^{l+1}=\frac{1}{M}$, for $i=1,\ldots,M$. Again, no forward model evaluations are required.
    \item \textsl{Mutation Step}: in order to improve the particle representation, we move the particles a bit with few applications of a $\mu_{\BY,l+1}$-invariant Markov chain Monte Carlo kernel. Each application requires a new forward model evaluation for each particle.
\end{enumerate}

We select the temperatures adaptively as described in \cite{jasra2011inference}, so as to preserve an effective sample size of $\frac{1}{1.1}M$ when selecting the new temperature. In the mutation step, we use random walk Metropolis-Hastings with a normally distributed increment in each dimension, with adaptively selected increment variances and number of Markov chain kernel applications \cite[Sect. 3.3]{beskos2015sequential}, the latter capped at $25$. Ideally, one would choose a Markov kernel robust with respect to the parameter dimension $J$, such as a preconditioned Crank-Nicolson (pCN) proposal in the Gaussian setting \cite{cotter2013mcmc,kantas2014sequential}. However, devising such dimension-robust proposals for non-Gaussian priors remains an open problem, with \cite{vollmer2015dimension} providing a first step for elliptic inverse problems, though its efficacy here is unclear, especially at high frequencies.

\subsection{Model and Algorithmic Settings}\label{ssec:settings}

In the following, we do not specify units for physical quantities as we assume we work with normalized ones. For our numerical experiments in two spatial dimensions, we consider the expansion \eqref{eq:radiusomega} with constant average radius $r_0\equiv 1$, Fourier basis
\begin{equation}\label{eq:fourierpsi}
    \psi_{2j-1}(\varphi)=\frac{1}{j}\sin(j\varphi),\quad\psi_{2j}(\varphi)=\frac{1}{j}\cos(j\varphi),
    \quad \varphi\in [0,2\pi),\, j\geq 1,
\end{equation}
and coefficients
\begin{equation}
    \beta_{2j-1}=\beta_{2j} = \frac{1}{4\sum^{\infty}_{k=1}\frac{1}{1+sk^{2+\epsilon}}}\cdot \frac{j}{\left(1+sj^{2+\epsilon}\right)},\quad j\geq 1,
    \label{eq:correlation}
\end{equation}
for a small $\epsilon>0$, that we took as $\epsilon=10^{-3}$. The parameter $s>0$ can be seen as a scaled version of the correlation length in the Whittle-Matérn covariance \cite[pp. 81--95]{williams2006gaussian}. In our numerical results, $s=0.2$. It is easy to check that the coefficients \eqref{eq:correlation} fulfill Assumptions \ref{ass:Expansion} and \ref{ass:uniform} with $\gamma_{\beta}=\frac{1}{2}$. Moreover, their decay ensures that $r\in C^1_{per}([0,2\pi))$ $\mu_0$-a.s. \cite[Lem. 2.5]{hiptmair2018large}.
We truncate the expansion \eqref{eq:radiusomega} to $J=10$ terms, such that $95\%$ of the variance of the infinite sum is covered.

In the forward problem \eqref{eq:HSTP}, we choose the case $\alpha_{in}=\alpha_{out}=1$. Together with our choice of the coefficient sequence \eqref{eq:correlation}, the first set of assumptions for Theorem \ref{thm:wellposed_HSTP} is fulfilled. The other material parameters are chosen as $n_{in}=0.5$, $n_{out}=1$, and the excitation is the plane wave $u^i(\Bx) = \exp\left(i\kappa_0\Bd\cdot\Bx\right)$ with $\Bd=\frac{1}{\sqrt{2}}(1,1)^{\top}$. We will test two values for the frequency of the incoming wave, $f_1 = 2\cdot 10^{8}$ and $f_2= 2f_1$, corresponding to the wavenumbers $\kappa_{0,1}=\frac{2\pi f_1}{c} = \frac{4}{3}\pi$ and $\kappa_{0,2}= \frac{8}{3}\pi$ respectively, with $c=3\cdot 10^{8}$ (m/s) the speed of light in vacuum.

The measurements are the real part of smoothed point values of the scattered field \eqref{eq:smoothedpointval} with $\sigma_o=0.25$ at $K=10$ points on a circle of radius $r_1=4$ around the scatterer,
\begin{equation}
	\mathbf{x}_i = \left(r_1\cos\left(\frac{2\pi i}{K}\right),r_1\sin\left(\frac{2\pi i}{ K}\right)\right)\quad\text{for }1\leq i \leq K.
	\label{eq:MeasuringPoints}
\end{equation}

To solve the forward problem on the reference configuration, we use an unstructured, simplicial mesh with maximum mesh size $h=0.15$ on the circular domain $B_R$ of radius $R=11$. For the lower frequency $f_1$, we use third order finite elements, for $f_2$ we will test fourth and fifth order finite elements (see Subsection \ref{ssec:numres} for more explanation). 

To solve the inverse problem, we consider measurements of the scattering behavior of a single object. We generate synthetic data by drawing a realization $\boldsymbol{y}^{\dagger}$ of $\boldsymbol{Y}^{\dagger}$ and solving the forward problem for the resulting ``true" shape with radius $r^{\dagger}$. We use the same realization for all numerical tests. The synthetic measurements have been computed using the same finite element order as the one used to sample the posterior, but on a staggered mesh with maximum mesh size $0.9 h$. The results have then been corrupted by noise, as specified in the next subsection.

In the sequential Monte Carlo algorithm, we use $M=1000$ particles. To plot the obtained marginal posteriors, we use a Gaussian kernel density estimation (Python function \texttt{scipy.gaussian\_kde} with default bandwidth selection).

\subsection{Numerical Results}\label{ssec:numres}

Our analytical results quantify the higher sensitivity to noise of the posterior at higher frequencies. To observe this, for each setting considered here we compare the posterior obtained from clean 
(noiseless) and noisy data. The estimates also involve wavenumber-independent constants affecting the posterior's sensitivity, as we now illustrate.

\subsubsection{The role of the Signal to Noise Ratio (SNR)}

We consider measurements of the scattered wave at the fixed locations \eqref{eq:MeasuringPoints}, which in the noisy case we corrupt by a realization of $\boldsymbol{\eta}\sim\mathcal{N}(\boldsymbol{0}, \sigma_{\eta}^2 \mathbb{I}_{10\times 10})$
with $\sigma_{\eta}^2 = 5\cdot 10^{-6}$ and $\mathbb{I}_{10\times 10}$ the identity matrix. We use the same noise realization for both frequencies. When running the inversion for the case $f_2= 2f_1$, we increase the order of the finite elements from three to four in order to keep a comparable level of discretization error in the forward problem \cite{bernkopf2025wavenumber}.

\begin{figure}
	\centering
	\begin{subfigure}{\textwidth}
	\centering
	\includegraphics[scale=0.6]{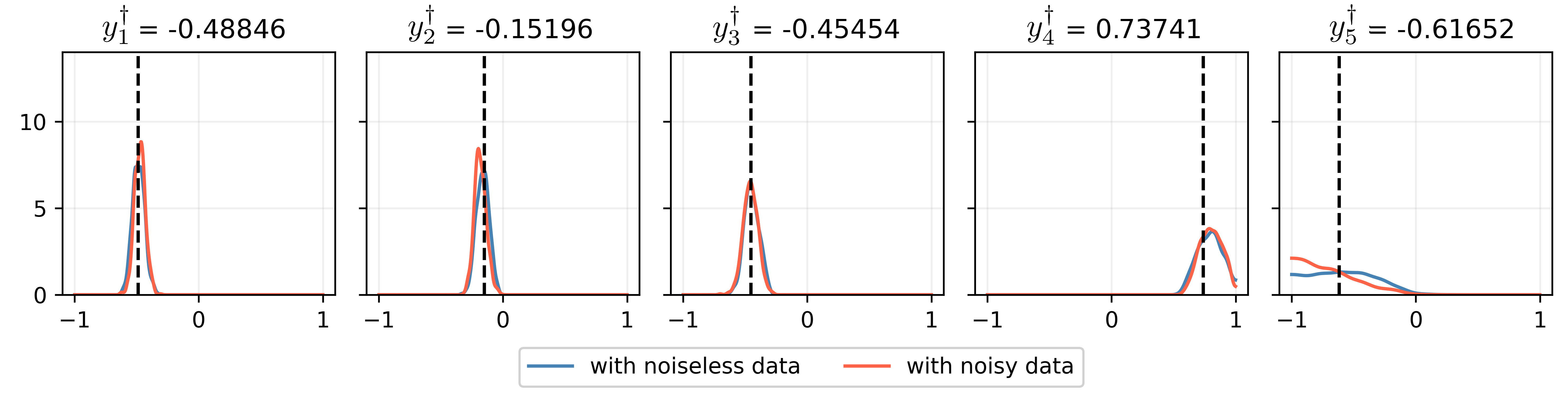}\caption{ }
\end{subfigure}%

\begin{subfigure}{\textwidth}
	\centering
	\includegraphics[scale=0.6]{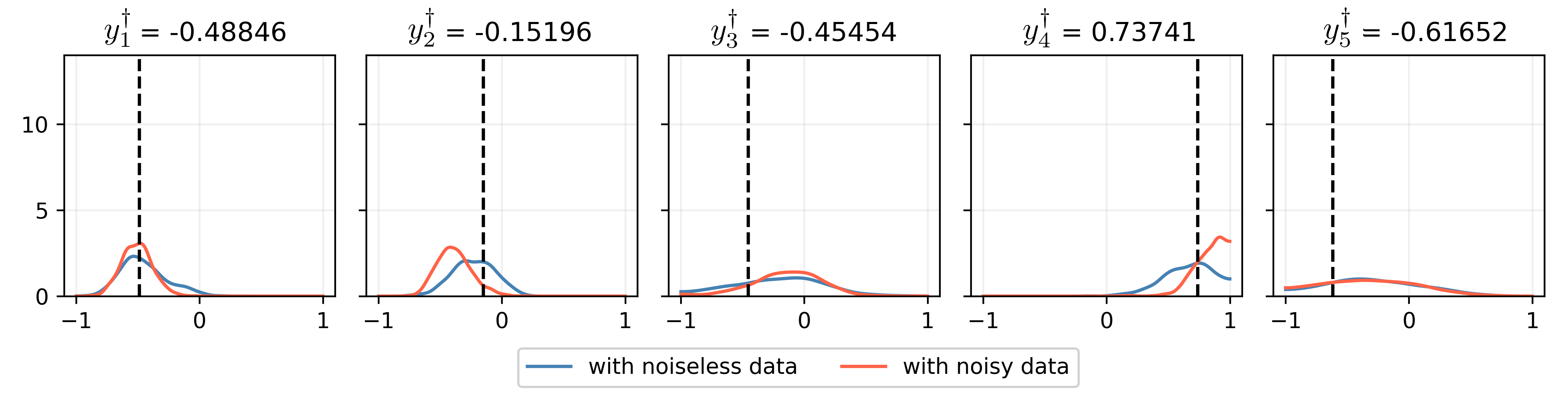}\caption{ }
\end{subfigure}
\caption{Marginal posteriors for the variables $\left\{Y_l\right\}_{l=1}^5$ in the expansion \eqref{eq:radiusomega}, using clean and noisy measurement data. (a): results for lower frequency, $f_1=2\cdot 10^{8}$, noise variance $\sigma_{\eta}^2=5\cdot 10^{-6}$ and using third order finite elements. (b): results for higher frequency, $f_2=2f_1$, noise variance $\sigma_{\eta}^2=5\cdot 10^{-6}$, the same noise realization as in (a) to generate the synthetic data, and fourth order finite elements. The dashed, vertical lines denote the values used to generate the synthetic data.}\label{fig:differentSNR_marginals}
\end{figure}

\begin{figure}
	\centering
	\begin{subfigure}[t]{0.5\textwidth}
		\centering
		\includegraphics[scale=0.75]{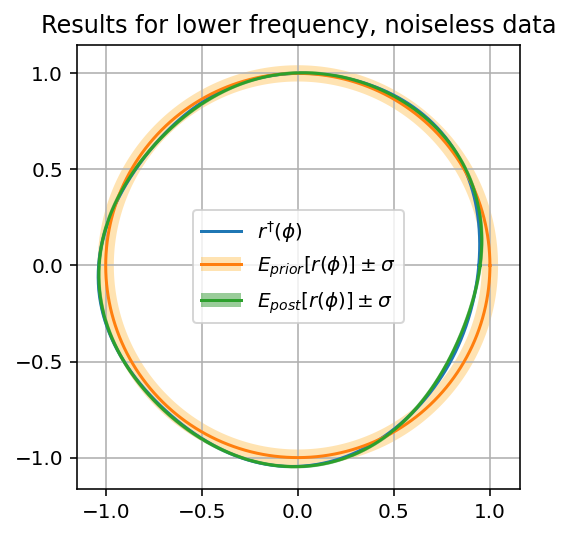}\caption{ }
	\end{subfigure}%
	~ 
	\begin{subfigure}[t]{0.5\textwidth}
		\centering
		\includegraphics[scale=0.75]{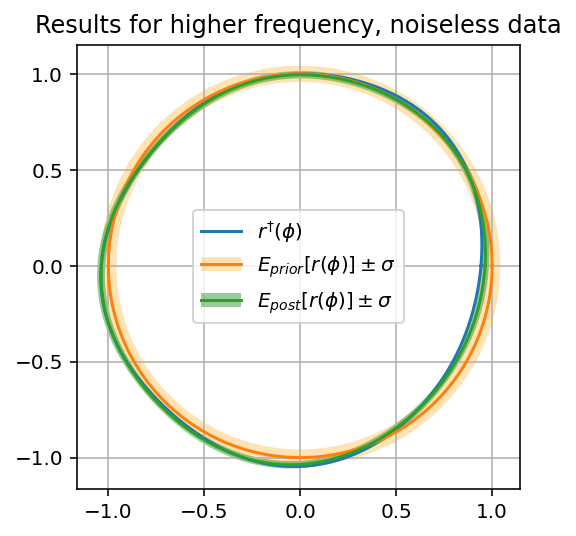}\caption{ }
	\end{subfigure}
	\caption{Posterior (and prior) averages $\pm$ standard deviation for the radius \eqref{eq:radiusomega} when using a plane wave excitation with frequency $f_1=2\cdot 10^{8}$ (a) and $f_2=2f_1$ (b), corresponding to the results in Figure \ref{fig:differentSNR_marginals} for noiseless data. Here $r^{\dagger}$ is the radius used to generate the data (the ``truth"). Noise realization and measurement locations are fixed for both frequencies, resulting in different SNR.}\label{fig:differentSNR_postavg}
\end{figure}

Figure \ref{fig:differentSNR_marginals} shows the marginal posteriors for the first five modes in the radius expansion \eqref{eq:radiusomega} and Figure \ref{fig:differentSNR_postavg} shows the posterior averages for the radius with standard deviations for the case of noiseless measurement data. We observe that the posteriors always concentrate around the true values used to generate the synthetic data. In Figure \ref{fig:differentSNR_postavg}(a), the uncertainty band around the posterior average (shaded green area) is so narrow that it is hardly visible. For the lower frequency, we do not observe a significant effect of the noise on the posterior: the blue and red densities in Figure \ref{fig:differentSNR_marginals}(a) are quite similar. For the higher frequency, we also do not observe much sensitivity of the posteriors to noise and the posteriors concentrate much less than in the case of lower frequency, as also slightly visible from Figure \ref{fig:differentSNR_postavg}(b). 

At first sight, these results seem contradicting the intuition and analytical results on higher sensitivity of the posterior to noise at higher frequencies. However, since the measurement locations and noise realization are fixed, and higher-frequency signals decay faster with distance from the scatterer, the measurements are smaller in magnitude, and thus less informative, at higher frequency: in our case, about one order of magnitude smaller for $f_2$ than for $f_1$.

Since our analytical estimates are not specialized for specific observation operators, this effect is not immediately evident from them. Still, if we inspect how they were obtained, e.g., from the proof of Corollary \ref{cor:pdewellposedness}, we see that what really matters in the sensitivity of the posterior to noise is the ratio between the smallest eigenvalue of the noise covariance ($\lambda_{min}$) and the noiseless data $\mathcal{G}(r)$, which corresponds essentially to the SNR ratio. 

\subsubsection{Sensitivity to the discretization error}

In this numerical test, we modify the noise realization and covariance matrix for $f_2$ so that the data have the same SNR as for $f_1$. 

Specifically, to compute the noise realization $\boldsymbol{\eta}_{2}$ for the higher frequency case and the associated noise covariance matrix, we take the noise realization and noise variance for the lower frequency case,  $\boldsymbol{\eta}_{1}$ and $\sigma_\eta=5\cdot 10^{-6}$, and rescale each entry of $\boldsymbol{\eta}_{1}$ and the diagonal entries of the (diagonal) noise covariance matrix as
\begin{equation}\label{eq:eta_adjusted}
	\eta_{2, i} = \eta_{1, i}\cdot\frac{\left|\mathcal{G}_{\kappa_{02},i}(r^{\dagger})\right|}{\left|\mathcal{G}_{\kappa_{01}, i}(r^{\dagger})\right|}, \qquad
	\Sigma_{i,i} = \sigma_{\eta} \cdot\frac{\left|\mathcal{G}_{\kappa_{02},i}(r^{\dagger})\right|}{\left|\mathcal{G}_{\kappa_{01}, i}(r^{\dagger})\right|},\quad i=1,\ldots, 10,
\end{equation}
where $\mathcal{G}_{\kappa_{01},i}(r^{\dagger})$ and $\mathcal{G}_{\kappa_{02},i}(r^{\dagger})$ denote the noiseless measurement at the $i$-th location when using an incoming wave with wavenumber $\kappa_{0,1} = 2\pi f_1$ and $\kappa_{0,2} = 2\pi f_2$, respectively.

\begin{figure}
	\centering
	\includegraphics[scale=0.6]{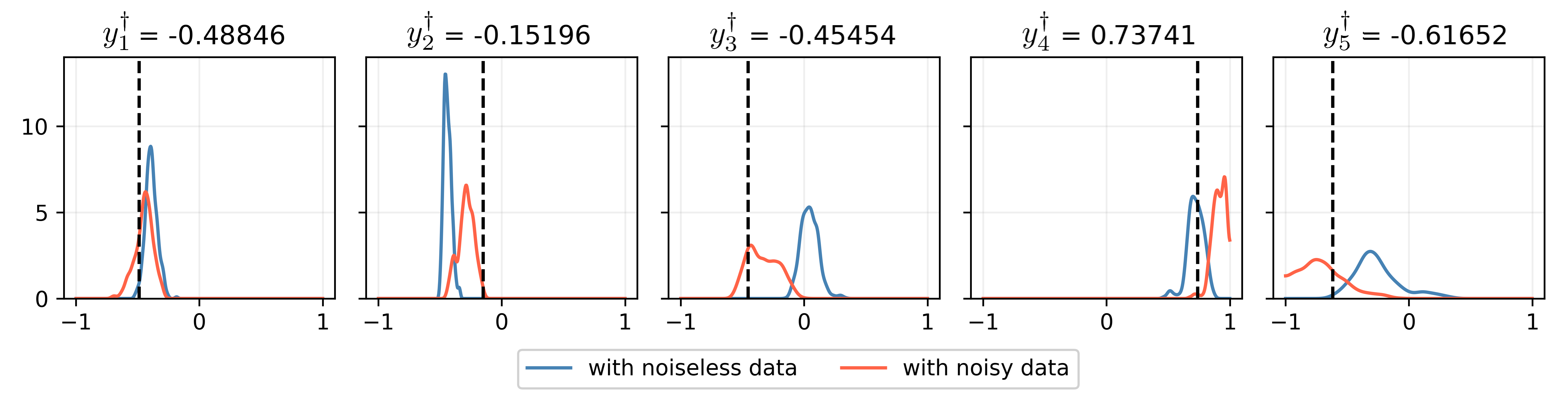}\caption{Marginal posteriors for the variables $\left\{Y_l\right\}_{l=1}^5$ in the expansion \eqref{eq:radiusomega}, using clean and noisy measurement data for the higher frequency case. The synthetically generated data without noise are the same as in Figure \ref{fig:differentSNR_marginals}, but now the noise realization and covariance have been adjusted according to \eqref{eq:eta_adjusted}, such that the SNR is the same here and in Figure \ref{fig:differentSNR_marginals}(a). }\label{fig:sameSNR_marginals}
\end{figure}

The first five marginal posteriors obtained when performing the inversion for the incoming wave with frequency $f_2$ are then depicted in Figure \ref{fig:sameSNR_marginals}. Not surprisingly, compared to Figure \ref{fig:differentSNR_marginals}(a), the posterior is more sensitive to noise. The interesting observation, however, is that the marginal posteriors obtained using noiseless data are often concentrated farther apart from the truth than those obtained from noisy data. 

This phenomenon can be attributed to the discretization error, despite the fact that we have used one order higher for finite elements in the higher frequency case to keep the discretization error in the forward model comparable to the lower frequency case. Since our analytical estimates are at the continuous level, the effect of the discretization error can then be interpreted as additional noise in the model output \cite{calvetti2018iterative}, see also discussion at the end of Subsection \ref{ssec:explbounds}. With the higher SNR used here, inspection of the generated data shows that the discretization error is of the same order as the added noise and, like the latter, has an effect on the posterior that is amplified at higher frequency, biasing the posterior obtained from noiseless data. Adjusting the finite element order according to the a priori forward-problem estimates in \cite{bernkopf2025wavenumber} is thus not sufficient. When the data is instead perturbed with noise, the noise can compensate for the discretization error and act as a regularizer, resulting in less biased posteriors. To validate this interpretation for the results in Figure \ref{fig:sameSNR_marginals}, we perform a last numerical test.

\subsubsection{Higher sensitivity to noise for same SNR and negligible discretization error}\label{sssec:adjustedSNRanddiscr}

In this test, we aim at comparing the posterior obtained using lower and higher frequency when: the SNR is the same in the two cases, and the discretization error is negligible. Since increasing the discretization accuracy to the point that, in the previous test, the effect of the discretization error is negligible compared to the statistical noise is very expensive computationally, we slightly increase the noise variance for the lower frequency, using $\sigma_{\eta}^2=2\cdot 10^{-5}$.

\begin{figure}
	\centering
	\begin{subfigure}{\textwidth}
		\centering
		\includegraphics[scale=0.6]{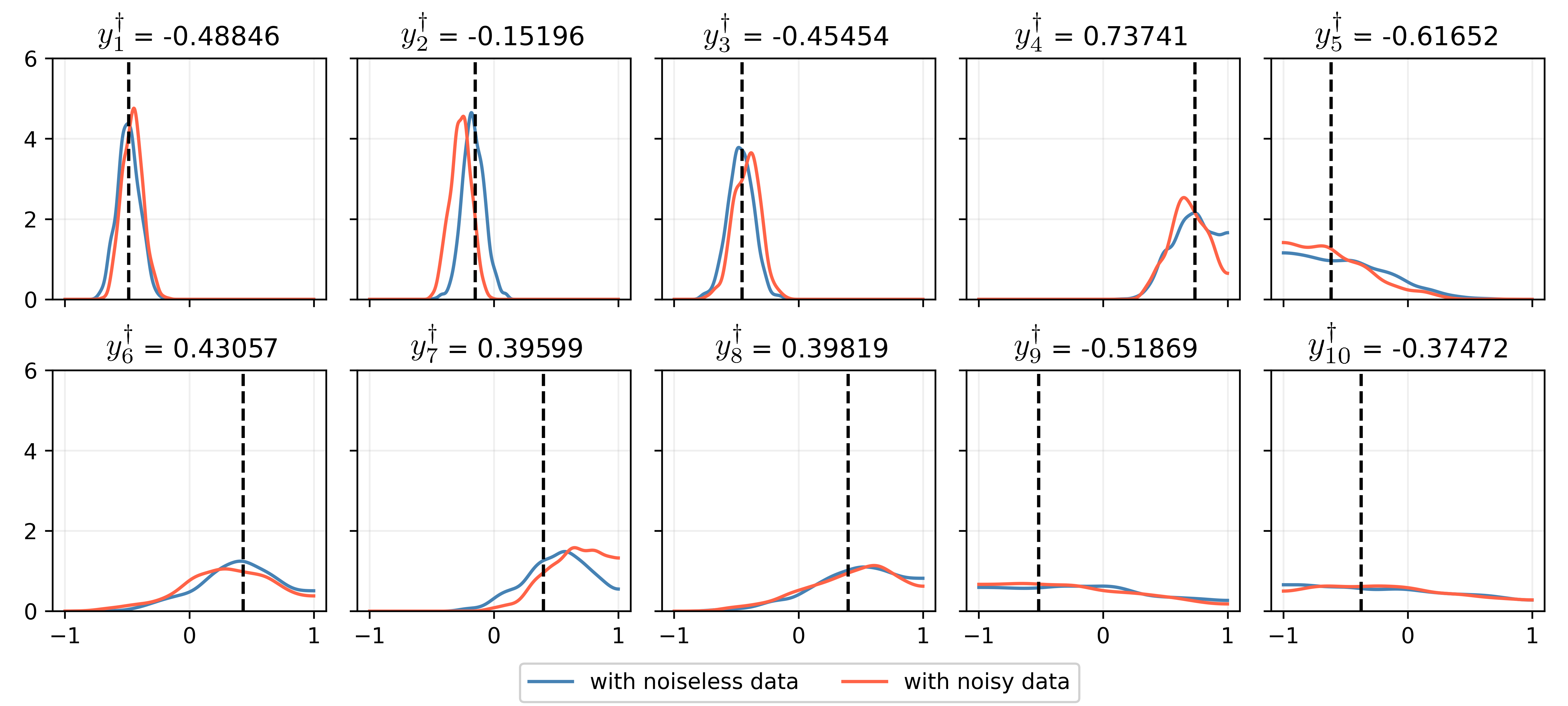}\caption{ }
	\end{subfigure}%
	
	\begin{subfigure}{\textwidth}
		\centering
		\includegraphics[scale=0.6]{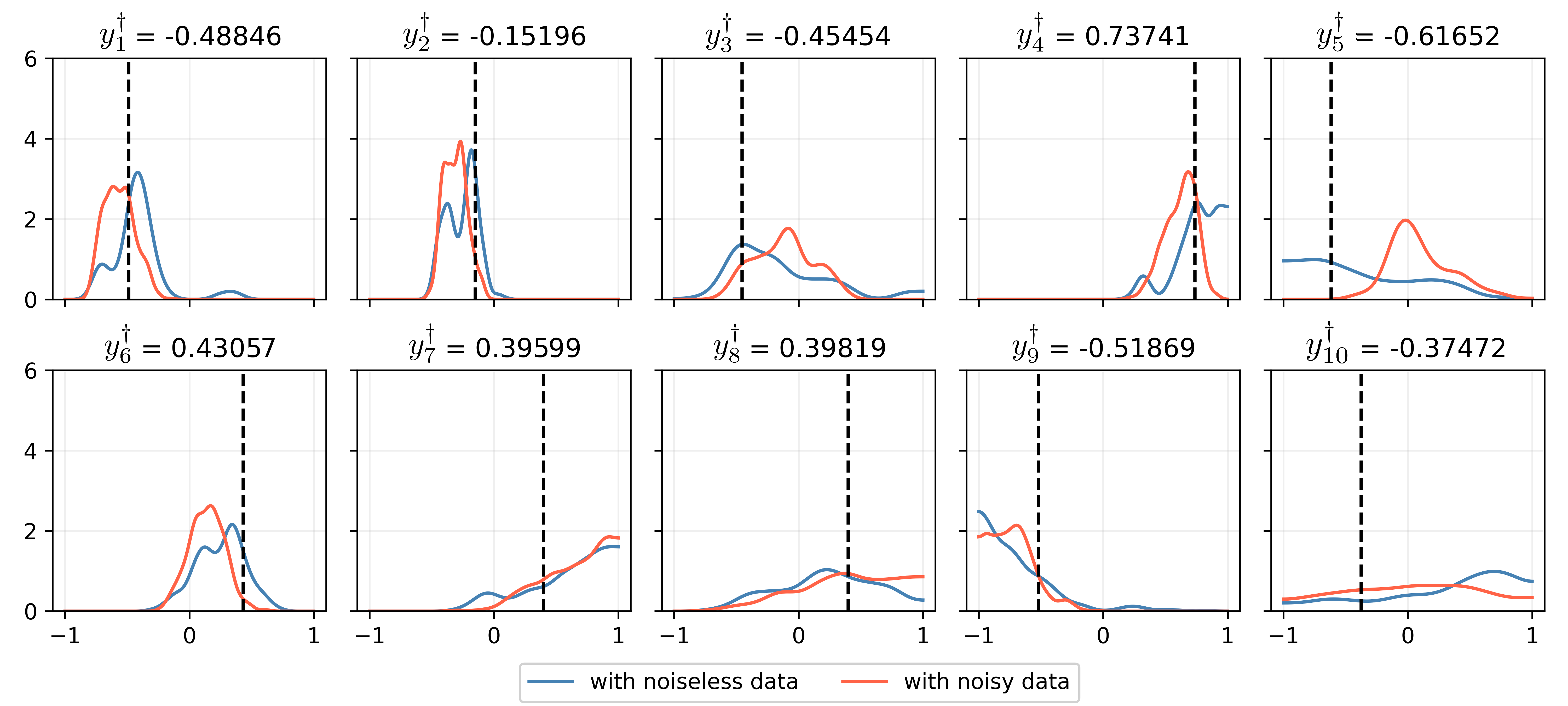}\caption{ }
	\end{subfigure}
	\caption{Marginal posteriors for the variables $\left\{Y_l\right\}_{l=1}^{10}$ in the expansion \eqref{eq:radiusomega}, using clean and noisy measurement data. (a): results for lower frequency, $f_1=2\cdot 10^{8}$, noise variance $\sigma_{\eta}^2=2\cdot 10^{-5}$ and using third order finite elements. (b): results for higher frequency, $f_2=2f_1$, noise realization and variance adjusted according to \eqref{eq:eta_adjusted}, and fifth order finite elements. The dashed, vertical lines denote the values used to generate the synthetic data.}\label{fig:sameSNRghigherdiscr_marginals}
\end{figure}

For the higher frequency case, we use fifth order finite elements (compared to fourth order used in the previous test), both to generate the data and for the inversion (we remind that the meshes are staggered to avoid the inverse crime); furthermore, we adjust the noise realization and covariance as in \eqref{eq:eta_adjusted}. The results are shown in Figure \ref{fig:sameSNRghigherdiscr_marginals}.

Now that the SNR is the same and that the effect of the discretization error on the posteriors is negligible in both cases, we see that, when the frequency is higher (Figure \ref{fig:sameSNRghigherdiscr_marginals}(b)), the posterior is more sensitive to noise. Furthermore, although the marginal posteriors for $Y_1, Y_2, Y_3$ and $Y_4$ seem to be slightly more concentrated in the lower frequency case, with the higher frequency we learn more on higher harmonics, mainly $Y_6$. Finally, the marginal posteriors in Figure \ref{fig:sameSNRghigherdiscr_marginals}(b) have more complex structure than in Figure \ref{fig:sameSNRghigherdiscr_marginals}(a), departing more from Gaussianity and exhibiting a slight multimodal behavior. Remarkably, all these effects are already visible at just twice the lower frequency, consistently with the sharp growth of the stability constants predicted by our estimates.

\section{Conclusions}\label{sec:conclusions}
While higher frequency in sensing technology is associated to higher precision, it also carries higher sensitivity to noise. This paper quantifies these effects in time-harmonic scattering problems modeled by the Helmholtz equation with unknown scatterer's shape. The main take-home messages are: 1) under some regularity assumptions, the Bayesian shape inverse problem is well-posed, and 2) sensitivity of the posterior to noise depends on the inverse relationship between the wavenumber and the geometrical scale of the problem (expressed for instance by terms like $\kappa_0 R$). The geometrical scale must also be taken into account when considering the measurements, as this may affect the signal to noise ratio as we have seen in the numerical results. Furthermore, the term noise should be interpreted in broader terms, including for instance a discretization or modeling error.

The results and techniques of this paper can be of interest to adjacent fields. They can be extended to the inverse problem for material parameters, e.g., for imaging purposes, starting from the estimates in \cite{graham2019helmholtz}, and the wavenumber-explicit bounds used in the proofs of Propositions \ref{prop:helmholtz_continuity1} and \ref{prop:soundsoft_continuity} can be relevant for shape optimization in acoustics and electromagnetics.

Several aspects need further investigation and will be addressed in further research, such as efficient numerical methods to sample the posterior at moderate to high frequencies and the treatment of functional data such as the far field.


\section*{Acknowledgements}
The authors thank E. Spence and A. Moiola for pointing them to the bound in Corollary \ref{cor:uscattbound} as improvement over the one in \cite{moiola2019acoustic}, and Wouter G. van Harten for assistance with the numerical implementation. LS thanks the organizers of the workshop \textsl{UQ in kinetic and transport equations and in high-frequency wave propagation} (2022 thematic program \textsl{Computational Uncertainty Quantification: Mathematical Foundations, Methodology \& Data}, ESI Vienna) for the scientific interactions that inspired this work, and Ralf Hiptmair for follow-up discussions.

LS was partially supported by NWO (Dutch Research Council), grant VI.Vidi.243.130.


\section*{Data Availability}
The code used to produce the numerical experiments is publicly available at \url{https://github.com/scalaura/Bayesianinversion_Helmholtzscattering/tree/main/Shape_inversion}.

\printbibliography

\end{document}